\documentclass[12pt]{article} 

\usepackage{euscript}
\usepackage{amsmath}
\usepackage{graphicx}
\usepackage[colors,economics]{optsys}
\usepackage{MyMacros}
\usepackage{stmaryrd}
\usepackage{tablefootnote}
\usepackage{dsfont}
\usepackage{caption}
\usepackage{subcaption}
\usepackage[dvipsnames]{xcolor}
\usepackage{eurosym} 
\usepackage{tikz}
\usetikzlibrary{decorations.pathmorphing}
\usetikzlibrary{decorations.pathreplacing,calligraphy}
\usepackage{algorithm2e}
\RestyleAlgo{ruled}
\SetKwComment{Comment}{/* }{ */}
\usepackage{soul}
\usepackage{hyperref}

\renewcommand{\WEEK}{\mathcal{S}}
\renewcommand{\HOUR}{\mathcal{H}}
\renewcommand{\STATE}{\mathcal{X}}

\renewcommand{\II}{\mathcal{I}}

\newcommand{\GRID}{G}

\usepackage{authblk}
\usepackage{amssymb}
\graphicspath{{Figures/}}
\newcommand{\mathscr}{\EuScript}
\usepackage{fullpage}

\newenvironment{keywords}{\bgroup
  \hsize=\textwidth%
  \noindent\unskip\textbf{Keywords.}\noindent\,\ignorespaces}
{\egroup}

\title{A Two-Timescale Decision-Hazard-Decision Formulation \\ for Storage Usage Values Calculation}

\author[1,2]{Camila Mart\'{\i}nez Parra}
\author[2]{Michel De Lara}
\author[2]{Jean-Philippe Chancelier}
\author[3]{Pierre Carpentier}
\author[1]{Jean-Marc Janin}
\author[1]{Manuel Ruiz}

\affil[1]{R\'eseau de Transport d'\'Electricit\'e, France}
\affil[2]{CERMICS, \'Ecole nationale des ponts et chauss\'{e}es, IP~Paris, France}
\affil[3]{UMA, ENSTA Paris, IP~Paris, France}

\date{\today}

\begin{document}

\maketitle

\begin{abstract}
  The penetration of renewable energies requires additional storages to deal with intermittency.
  Accordingly, there is growing interest in evaluating the usage value (opportunity cost) associated
  with stored energy in large storages, a cost
  obtained by solving a multistage stochastic optimization problem.
  Today, to compute usage values under uncertainties,
  an adequacy resource problem is solved using stochastic dynamic programming assuming a
  {hazard-decision} information structure.
  This modelling assumes complete knowledge of the coming week uncertainties, which
  is not adapted to the system operation as they occur at smaller timescale (hour). 
  This is why, in this paper, we propose to model the decision-making process
  with, on the one hand, two timescales (weekly time decomposition, but
  uncertainties that are realized hour by hour, and hourly energy balance
  constraints) and, on the other hand,
  a {decision-hazard-decision} information structure considering both planning
  and recourse  decisions.
  This structure is used to decompose the multistage  decision-making process
  into a nonanticipative  planning step in which  the binary on/off decisions for
  the thermal units are made, and a recourse step in which the power modulation decisions
  are made once the uncertainties have been disclosed.
  We perform time decomposition by stochastic dynamic programming, and not by
  SDDP, due to the presence of binary on-off decision variables.
  In a numerical case, we illustrate how usage values are sensitive to how
  the disclosure of information is modelled.
\end{abstract}



\section{Introduction}

In~\S\ref{Context:_usage_values_in_prospective_studies},
we outline the context of usage values in prospective studies.
In~\S\ref{Our_contribution}, we sketch the main features of our contribution and then,
in~\S\ref{Comparison_with_the_literature}, we make comparisons with the literature.
We detail the structure of the paper in~\S\ref{Structure_of_the_paper}.


 \subsection{Context: usage values in prospective studies}
\label{Context:_usage_values_in_prospective_studies}

In energy systems, prospective studies aim at identifying possible curtailment in demand or production,
network congestions, or the non-satisfaction of greenhouse gas emission targets.
This is done by means of simulations.
The simulation of an energy system operation is carried out by solving its associated \emph{adequacy resource} problem.
Mathematically, the adequacy resource problem is a \emph{multistage optimization} problem that aims to allocate,
 hour by hour, the dispatchable production means such that the demand is met while minimizing the overall production cost.
 The problem  is formulated from the point of view of a central planner who takes all decisions in the system
 for the common good.
 In such case, the production costs are fixed and resources are allocated by merit order, i.e.
  the cheapest units are used before the most expensive ones.
  Since the electrical demand and the availability of the thermal units are uncertain, the adequacy resource problem is
\emph{stochastic}.
This problem is naturally formulated using two timescales where the planning is done  in coordinated manner for the week
ahead but the hourly energy balance has to be met, that is, a weekly timescale and an hourly timescale.
In this context, the question of when the energy in storages is going to be used arises.
\emph{Usage values} are the storages' prices --- a price signal that makes it
possible to choose when and how much of its energy is used --- that depend on the energy
system setting.

We now discuss how usage values are computed for a system composed of several storages and thermal units.
Stochasticity is introduced by the residual demand of the system
(difference between demand and non-dispatchable production)
and the availability of the thermal units (outages).
We do not focus on how the uncertainty scenarios are generated, but on \emph{information structure}, that is,
on how the available information about these scenarios is disclosed throughout the decision-making process.
The dynamics in the storages introduce temporal coupling in the problem, leading to a
\emph{multistage stochastic optimization problem}. A standard way to solve it
is by dynamic programming. Then, we obtain usage values by differentiating the Bellman value functions with respect to
the storage levels.
%
%
%
%
Once the usage values are calculated, one can carry out prospective studies by simulation.
Indeed, one can compute a storage's management policy to be used in the
resolution of the adequacy resource problem.
However, it is important to highlight that the adequacy resource problem studied here
does not intend to yield an implementable operational schedule for the system,
but to give an overview of the system operation.

In practice, the French TSO R\'{e}seau de Transport d'Electricit\'{e} (RTE)\footnote{https://www.rte-france.com/}
is dealing with such prospective studies using the open-source tool Antares~\cite{Antares}.
RTE tackles the resolution of the {multistage stochastic optimization problem} associated with the adequacy problem
as follows. The timespan is one year, with a weekly timescale and an hourly timescale.
Dynamic programming is performed at the weekly timescale: a Bellman value function is computed at the beginning of each week,
anticipating the coming uncertainties over the week and respecting the hourly energy balance equations within the week.
This structure assumes that all the decisions in a week are made with full knowledge of the uncertainties of the week.
When the dispatchable units are ``fast" to start (like hydropower or gas turbine),
it is reasonable to assume that we can ``wait-and-see" the uncertainties
before making the decisions.
By contrast, as ``slow" dispatchable units (like nuclear or coal) need more time to start producing,
start or stop decisions must be made before the uncertainties are known (especially the outages),
and ``wait-and-see" decisions are not adapted.
This is why, our main contribution focuses on informations structures suited to handle both
fast and slow dispatchable units submitted to outages,
with two timescales.

\subsection{Our contributions}
\label{Our_contribution}

Information structures are used to model the information available at each stage of the decision-making process.
The current practice at RTE, as described just above, corresponds to what we call \emph{weekly hazard-decision}.
Indeed, this structure assumes that all the decisions in a week are made with full knowledge of the uncertainties of the week.

In~\S\ref{Comparison_with_the_literature}, we present papers that deal with so-called
\emph{hazard-decision}, \emph{decision-hazard} and \emph{full decision-hazard}, especially in relation to the
SDDP algorithm \cite{Pereira-Pinto:1991}.
By comparison, our contributions in this paper are the following.
\begin{itemize}
\item
  We study \emph{decision-hazard-decision} information structures \emph{with two timescales}
  (already introduced in \cite{Carpentier-Chancelier-DeLara-Martin-Rigaut:2023}),
to address the usage value calculation by keeping track of hourly constraints but allowing a
  weekly decomposition of the yearly problem.
\item
  We consider a model with binary variables to model the on/off decisions for the thermal (dispatchable) units,
  and solve the stochastic dynamic programming (SDP) equations in a decision-hazard-decision framework without relaxing these
  binary decisions.
\end{itemize}

Regarding this second item, the energy systems that we consider
include, in addition to storage, thermal units that may be slow to go
from off to on (and reverse), and which may be subject to uncertain
outages.  The on/off decisions are naturally modelled as binary
decisions and, since we consider a problem with a small space state
dimension, we choose to use SDP (with special versions of Bellman
equations for nonanticipative information structures at the week
timescale).  The choice of using SDP is on purpose, and lies in the
fact that we seek to present a framework capable of handling binary or
discrete variables within the planning decisions (switch on/off), in
which case SDDP is not suitable, although variants like SDDiP \cite{zou2019stochastic}
exist. The scalability of the problem concerning the
dimension of the state space will be addressed in future work by
applying spatial decomposition techniques.  \medskip

We now sketch what we mean by {decision-hazard-decision} information
structures {with two timescales}.  First, we describe the ideal case
of an information structure with weekly planning decisions and hourly recourses, as illustrated in
Fig.~\ref{Sketch_of_information_structure_with_weekly_planning_decisions_and_hourly_recourses}.
We do not consider that, when we make a (planning) decision at the
beginning of the week, we know in advance the hourly random variables
that will materialize during the coming week.  Rather, at the end of
every hour, we allow for an hourly decision variable (recourse) to
handle hourly balance constraints, knowing all uncertainties up to
this hour.

  \begin{figure}[htbp]
    \centering
    \resizebox{0.8\columnwidth}{!}{
      \begin{tikzpicture}
       \draw[dashed, thick] (-6.5,0) -- (-5.5,0);
       \draw [ -latex,  dashed] (-5.5,1)-- (-5.5,0.15);
       \node (us) at (-5.5,1.4) {\footnotesize {$\substack{\mathsf{D}\text{ weekly}\\\text{(planning)}}$}};
       \filldraw[black] (-5.5,0) circle (2pt) node[anchor=north] {$0$};
       \draw[black, thick] (-5.5,0) -- (-4,0);
       \filldraw[black] (-4,0) circle (2pt) node[anchor=north] {$1$};
       \draw[dashed, thick] (-4,0) -- (-2.5,0);
       \filldraw[black] (-2.5,0) circle (2pt) node[anchor=north] {$\hour-1$};
       \draw[black, thick] (-2.5,0) -- (-1,0);
       \filldraw[black] (-1,0) circle (2pt) node[anchor=north] {$\hour$};
       \draw[black, thick] (-1,0) -- (0.5,0);
       \filldraw[black] (0.5,0) circle (2pt) node[anchor=north] {$\hour+1$};
       \draw[dashed, thick] (0.5,0) -- (2,0);
       \filldraw[black] (2,0) circle (2pt) node[anchor=north] {$167$};
       \draw[dashed, thick] (0.5,0) -- (2,0);
       \draw[black, thick] (02 ,0) -- (3.5,0);
       \filldraw[black] (3.5,0) circle (2pt) node[anchor=north] {$168$};
       \draw[dashed, thick] (3.5,0) -- (4.5,0);
            \draw[ -latex,  dashed] (-4 ,1)-- (-4, 0.15);
            \node (1) at (-4,1.2) {\footnotesize {$\substack{\mathsf{D}\text{ hour 1}\\\text{(recourse)}}$}};
            \draw[ -latex,  dashed] (-2.5 ,1)-- (-2.5,0.15);
            \node (1) at (-2.5,1.2) {\footnotesize {$\substack{\mathsf{D}\text{ hour $\hour$-1}\\\text{(recourse)}}$}};
            \draw[ -latex,  dashed] (-1 ,1)-- (-1,0.15);
            \node (1) at (-1,1.2) {\footnotesize {$\substack{\mathsf{D}\text{ hour $\hour$}\\\text{(recourse)}}$}};
            \draw[ -latex,  dashed] (0.5 ,1)-- (0.5,0.15);
            \node (1) at (0.5,1.2) {\footnotesize {$\substack{\mathsf{D}\text{ hour $\hour$+1}\\\text{(recourse)}}$}};
            \draw[ -latex,  dashed] (2 ,1)-- (2,0.15);
            \node (1) at (2,1.2) {\footnotesize {$\substack{\mathsf{D}\text{ hour $167$}\\\text{(recourse)}}$}};
            \draw[ -latex,  dashed] (3.5 ,1)-- (3.5,0.15);
            \node (1) at (3.5,1.2) {\footnotesize {$\substack{\mathsf{D}\text{ hour $168$}\\\text{(recourse)}}$}};
       \draw [ -latex,  dashed] (-5.5,1)-- (-4,0.15);
       \draw [ -latex,  dashed] (-5.5,1)-- (-2.5,0.15);
       \draw [ -latex,  dashed] (-5.5,1)-- (-1,0.15);
       \draw [ -latex,  dashed] (-5.5,1)-- (0.5,0.15);
       \draw [ -latex,  dashed] (-5.5,1)-- (2,0.15);
  \end{tikzpicture}
}
\caption{Sketch of information structure with weekly planning decisions and
  hourly recourses.
An arrow maps the available information towards the decision,
  so that, here, the information structure is nonanticipative as all arrows go either down or from the left to the right.
\label{Sketch_of_information_structure_with_weekly_planning_decisions_and_hourly_recourses}}
\end{figure}

We denote this information structure with weekly planning decisions
and hourly recourses by~$\WPHR$ (weekly~D and hourly~HD), and
symbolically represent it
by~\eqref{eq:Sketch_of_information_structure_with_weekly_planning_decisions_and_hourly_recourses}:
\begin{equation}
  \WPHR :
  \overbrace{\mathsf{D}-\underbrace{\mathsf{HD}-\cdots-\mathsf{HD}}_{\text{168
        hours}}}^{\text{one week}}
  \eqfinp
  \label{eq:Sketch_of_information_structure_with_weekly_planning_decisions_and_hourly_recourses}
\end{equation}
The (backward) Bellman equation 
corresponding to~$\WPHR$ can be established
from~\cite{Carpentier-Chancelier-DeLara-Martin-Rigaut:2023}, and is sketched
in~\eqref{eq:Sketch_of_information_structure_with_weekly_planning_decisions_and_hourly_recourses_Belmann}
(see also \S\ref{Bellman_equations_in_decision-hazard-decision_information_structure_with_hourly_recourse}):
{\small 
\begin{align}
  & \substack{\text{cost}\\\text{to go}\\ \WPHR } =
  \label{eq:Sketch_of_information_structure_with_weekly_planning_decisions_and_hourly_recourses_Belmann}
  \\
  &
    \min_{\substack{ \mathsf{D} \text{ weekly}\\\text{(planning)}} }
  \EE  \Biggl[
  {         \min_{\substack{\mathsf{D} \text{ hour 1}\\\text{(recourse)}} }}
  \EE  \biggl[
  \dots
  {               \min_{\substack{\mathsf{D} \text{ hour 167}\\\text{(recourse)}} } }
  \EE\Bigl[
  {           \min_{\substack{\mathsf{D} \text{ hour 168}\\\text{(recourse)}} } }
  \bigl(
  \substack{\text{weekly}\\\text{cost}}  
  +
  \substack{\text{ next cost}\\\text{to go}\\ \WPHR }
  \bigr)
  \mid \mathsf{H}_1,  \mathsf{H}_2 ,\dots,  \mathsf{H}_{167}
  \Bigr]
  \dots
  \mid   \mathsf{H}_1\biggr]
  \Biggr]
  \nonumber
  \eqfinp
\end{align}
}
The obtained Bellman functions are optimal
under a (debatable) assumption of statistical independence of the random
variables (weather conditions, demand, outages, etc.) between weeks.
As commonly done with stochastic dual dynamic programming (SDDP),
we could account for lag variables to be added aside storage,
hence making an extended state.
Of course, as we use SDP, we cannot extend the one dimensional state (storage level)
with more than three or four entries.
In the paper, we choose to keep only the storage as state.

Even without possible lag variables,
the Bellman equation sketched
in~\eqref{eq:Sketch_of_information_structure_with_weekly_planning_decisions_and_hourly_recourses_Belmann}
is numerically untractable for several reasons:
because of the DH part in the information structure, we cannot use SDDP without adaptation;
due to the presence of binary decision variables (on-off), we cannot resort to
the SDDP algorithm;
as we do not suppose statistical independence of the hourly random variables
inside a week, the intra-week problem is out of numerical reach.

This is why our approach consists in relaxing the problem by collecting the hourly decision variables
inside the week in a single vector decision (recourse) made at the end of the
week. With this, the DH phase is now followed by a
hazard-decision (HD) phase, making the whole information structure
\emph{decision-hazard-decision at the weekly timescale} ($\DHD$), as illustrated in
Fig.~\ref{Decision-hazard-decision_information_structure_at_the_weekly_timescale}.
\begin{figure}[htbp]
  \centering
  \resizebox{0.8\columnwidth}{!}{
    \begin{tikzpicture}
     \draw[dashed, thick] (-6.5,0) -- (-5.5,0);
     \draw [ -latex,  dashed] (-5.5,1)-- (-5.5,0.15);
     \node (us) at (-5.5,1.4) {\footnotesize {$\substack{\mathsf{D}\text{ weekly}\\\text{(planning)}}$}};
     \filldraw[black] (-5.5,0) circle (2pt) node[anchor=north] {$0$};
     \draw[black, thick] (-5.5,0) -- (-4,0);
     \filldraw[black] (-4,0) circle (2pt) node[anchor=north] {$1$};
     \draw[dashed, thick] (-4,0) -- (-2.5,0);
     \filldraw[black] (-2.5,0) circle (2pt) node[anchor=north] {$\hour-1$};
     \draw[black, thick] (-2.5,0) -- (-1,0);
     \filldraw[black] (-1,0) circle (2pt) node[anchor=north] {$\hour$};
     \draw[black, thick] (-1,0) -- (0.5,0);
     \filldraw[black] (0.5,0) circle (2pt) node[anchor=north] {$\hour+1$};
     \draw[dashed, thick] (0.5,0) -- (2,0);
     \filldraw[black] (2,0) circle (2pt) node[anchor=north] {$167$};
     \draw[dashed, thick] (0.5,0) -- (2,0);
     \draw[black, thick] (02 ,0) -- (3.5,0);
     \filldraw[black] (3.5,0) circle (2pt) node[anchor=north] {$168$};
     \draw[dashed, thick] (3.5,0) -- (4.5,0);
    \draw[ -latex,  dashed] (3.5 ,1)-- (-4, 0.15);
    \draw[ -latex,  dashed] (3.5 ,1)-- (-2.5,0.15);
    \draw[ -latex,  dashed] (3.5 ,1)-- (-1,0.15);
    \draw[ -latex,  dashed] (3.5 ,1)-- (0.5,0.15);
    \draw[ -latex,  dashed] (3.5 ,1)-- (3.5,0.15);
    \node (1) at (03.5,1.4) {\footnotesize {$\substack{\mathsf{D}\text{ weekly}\\\text{(recourse)}}$}};
    \draw [ -latex,  dashed] (-5.5,1)-- (-4,0.15);
    \draw [ -latex,  dashed] (-5.5,1)-- (-2.5,0.15);
    \draw [ -latex,  dashed] (-5.5,1)-- (-1,0.15);
    \draw [ -latex,  dashed] (-5.5,1)-- (0.5,0.15);
    \draw [ -latex,  dashed] (-5.5,1)-- (2,0.15);
\end{tikzpicture}
}
\caption{Decision-hazard-decision information structure at the weekly timescale.
  An arrow maps the available information towards the decision,
  so that, here, the information structure is partly anticipative as some arrows go from the right
  to the left.
\label{Decision-hazard-decision_information_structure_at_the_weekly_timescale}}
\end{figure}

This relaxed version condenses all the hourly recourse decisions in one vector decision, leading
to a more manageable expression for the Bellman equations sketched
in~\eqref{Decision-hazard-decision_information_structure_at_the_weekly_timescale_Bellman}
(see also~\S\ref{subsection:BellmanDHD}):
\begin{align}
  &
  \substack{\text{cost}\\\text{to go} \\ \DHD}
  =   \min_{\substack{ \mathsf{D} \text{ weekly}\\\text{(planning)}} }
               \EE  \Biggl[
                {         \min_{\substack{\mathsf{D} \text{ weekly}\\\text{(recourse)}} }}
           \biggl(
                        \substack{\text{weekly}\\\text{cost}}  
                           +
                           \substack{\text{next cost}\\\text{to go} \\ \DHD}
             \biggr)
  \Biggr]
 \eqfinp
  \label{Decision-hazard-decision_information_structure_at_the_weekly_timescale_Bellman}
\end{align}
The DH phase (weekly planning) is now followed by a
hazard-decision (HD) phase (weekly recourse), making the whole information structure
decision-hazard-decision.
With this approach, we are able to handle a $\DHD$ information structure
without having to extend the state by adding a control aside,
contrarily to other approaches (see discussion in \S\ref{Street-Valladao-Lawson-Velloso:2020}).

\subsection{Comparison with the literature}
\label{Comparison_with_the_literature}
In~\S\ref{Street-Valladao-Lawson-Velloso:2020},
~\S\ref{Valladao-Silva-Poggi:2019}
and~\S\ref{dowson2020policy}, we sketch the contribution of three papers that deal with decision-hazard
and hazard-decision information structures in multistage stochastic optimization problems.
We present them par decreasing order of proximity with our work.
Finally, in~\S\ref{Synthesis_of_our_contribution_in_comparison_with_the_literature},
we posit our contribution in comparison with the literature in a synthetic fashion.

\subsubsection{%
A.~Street, D.~Vallad\~{a}o, A.~Lawson, and A.~Velloso.
\cite{Street-Valladao-Lawson-Velloso:2020} {\em Applied Energy}, 2020}
\label{Street-Valladao-Lawson-Velloso:2020}

The paper \cite{Street-Valladao-Lawson-Velloso:2020} is a nice one whose title \emph{Assessing the cost of the hazard-decision simplification in
  multistage stochastic hydrothermal scheduling} clearly indicates the interest of the authors to discuss the impact of
the hazard-decision simplification in multistage stochastic optimization problems.
We share the same preoccupation as \cite{Street-Valladao-Lawson-Velloso:2020}.

In \cite[Sect.~3]{Street-Valladao-Lawson-Velloso:2020}, the authors focus on the handling of DH
information structure in multistage stochastic optimization problems
formulated with continuous decision variables,
convex costs and constraints, and linear dynamics.
Their main preoccupation is how the SDDP algorithm can be adapted to the DH setting.
Thus, they make a stagewise independence assumption
and state Bellman equations without hint at a proof.

Starting from the Bellman equations \cite[Equations~(6)-(11)]{Street-Valladao-Lawson-Velloso:2020} in HD, they propose two
reformulations to handle DH.
Decision variables are split in two subgroups:
those made under uncertainty (also called here-and-now, preventive)
and those made after the observation (wait-and-see, corrective).
\begin{itemize}
\item
  The fullDH-SDDP Equations \cite[Equations~(12)-(17)]{Street-Valladao-Lawson-Velloso:2020} --- although called DH formulation ---
  can be interpreted as Bellman equations in DHD.
  Indeed, in the splitting of the decision variables in two subgroups,
  the decisions made under uncertainty (here-and-now, preventive)
  are the first D in the DHD information structure, hence
  are subject to the nonanticipativity constraints~\cite[Equation~(16)]{Street-Valladao-Lawson-Velloso:2020}
  in a two stage stochastic optimization problem.
  The decisions made after the observation (wait-and-see, corrective) are the last
D in the DHD information structure, hence correspond to recourse decision variables
in a two stage stochastic optimization problem.
  \item
    ASDH-SDDP Equations \cite[Equations~(18)-(23)]{Street-Valladao-Lawson-Velloso:2020}
    and \cite[Equations~(24)-(31)]{Street-Valladao-Lawson-Velloso:2020} are Bellman equations for value functions
    whose argument is an extended state (nonanticipative controls are added
    aside the original state to make an extended state).
\end{itemize}
Our approach differs from ASDH-SDDP in that we do not extend the state.
Our approach is close to fullDH-SDDP, but is not restricted to
continuous decision variables, convex costs and constraints, and linear
dynamics.
We provide Bellman equations --- whose derivation is not given in the paper, as
it is lengthy, but that can be established from~\cite{Carpentier-Chancelier-DeLara-Martin-Rigaut:2023} ---
for general decision variables (continuous, discrete or mixed),
costs and constraints (convex or not), and dynamics (linear or not).
In~\cite{Street-Valladao-Lawson-Velloso:2020}, Bellman equations are stated in the SDDP setting,
that is, for continuous decision variables, convex costs and constraints, and linear
dynamics.
Moreover, we give a formal definition of DHD, independently of SDDP,
whereas \cite{Street-Valladao-Lawson-Velloso:2020} evokes DH and HD in the SDDP setting.
Finally, we tackle the issue of DHD with two timescales, whereas
  \cite{Street-Valladao-Lawson-Velloso:2020} deals with a single timescale.

Regarding the energy application, the authors focus on operation --- that is, on finding strategies to operate
hydrothermal scheduling --- whereas we focus on usage values in prospective studies.

This said, let us insist that \cite{Street-Valladao-Lawson-Velloso:2020} is a nice contribution in raising awareness about the cost of the hazard-decision simplification in
multistage stochastic optimization problems in energy systems --- where
slow activation thermoelectric units need to be scheduled as here-and-now decisions,
while fast-activating thermoelectric and hydroelectric units are considered as wait-and-see decisions.

\subsubsection{%
D.~Vallad\~{a}o, T.~Silva, and M.~Poggi.
\cite{Valladao-Silva-Poggi:2019}
{\em Annals of Operations Research}, 2019}
\label{Valladao-Silva-Poggi:2019}

To the difference with~\cite{Street-Valladao-Lawson-Velloso:2020}, the paper
\cite{Valladao-Silva-Poggi:2019} is focused neither on energy systems, nor on the DH information structure as no more than
two pages are devoted to this issue.
However, in \cite[Sect.~4]{Valladao-Silva-Poggi:2019}, the authors consider a multistage
  stochastic optimization problem with risk constraints at every stage and,
  because of such a constraint, it would not make sense to know the uncertainty
  in advance, which justifies a DH information structure.
  As the portfolio problem under study is formulated with continuous decision variables,
  convex costs and constraints, and linear dynamics, the authors resort to SDDP as resolution method.
  Then, the authors discuss how to adapt SDDP to a DH formulation for the
  Markov chained stochastic dual dynamic programming (MSDDP) algorithm.
  They propose to solve a two-stage problem for each SDDP subproblem in
  the risk-constrained dynamic portfolio optimization context.
  A Bellman equation is given in the linear case, without hint at a proof.

  Thus, one year before~\cite{Street-Valladao-Lawson-Velloso:2020} was published, \cite{Valladao-Silva-Poggi:2019} raises
  the point that the DH information structure should be considered,
  a point with which we fully agree. Then, \cite{Valladao-Silva-Poggi:2019} proposes a way to
  adapt the SDDP algorithm, without state extension (hence keeping the original state as argument of Bellman functions).
  The approach is that of fullDH-SDDP in~\cite{Street-Valladao-Lawson-Velloso:2020}.

  The comparison of \cite{Valladao-Silva-Poggi:2019} with our work is much more limited than with~\cite{Street-Valladao-Lawson-Velloso:2020},
  since \cite{Valladao-Silva-Poggi:2019} is not really focused on the DH information structure.
  Contrarily to~\cite{Valladao-Silva-Poggi:2019}, our approach is not restricted to
continuous decision variables, convex costs and constraints, and linear dynamics.
We provide Bellman equations --- whose derivation is not given in the paper, as
it is lengthy, but that can be established from~\cite{Carpentier-Chancelier-DeLara-Martin-Rigaut:2023} ---
for general decision variables (continuous, discrete or mixed),
costs and constraints (convex or not), and dynamics (linear or not).
In~\cite{Valladao-Silva-Poggi:2019}, Bellman equations are stated in the SDDP setting,
that is, for continuous decision variables, convex costs and constraints, and linear
dynamics.
Finally, we tackle the issue of DHD with two timescales, whereas
  \cite{Valladao-Silva-Poggi:2019} deals with a single timescale.

\subsubsection{%
 O.~Dowson.
\cite{dowson2020policy}
 {\em Networks}, 2020}
\label{dowson2020policy}

To the difference with~\cite{Street-Valladao-Lawson-Velloso:2020}, the paper~\cite{dowson2020policy}
is not focused on the DH information structure as it is only tackled in slighly more than one
page in~\cite[\S2.2 (Nodes)]{dowson2020policy}, and then in half page in~\cite[\S3.4 (Standart form)]{dowson2020policy}.
In~\cite[\S2.2]{dowson2020policy}, the author adopts
the formalism of stochastic optimal control in discrete time.
This done, HD and DH are defined in the framework of state feedback policies:
HD is when a policy depends both on current state and uncertainty,
whereas DH is when a policy depends only on the current state.
Our approach is more general in that the DHD structure is expressed in terms of measurability
constraints, using \sigmafields generated by the uncertainty process.
At the end of~\cite[\S2.2]{dowson2020policy}, the author indicates that it is possible to tranform a DH node into a
deterministic node followed by a HD one with an extended state space
(as in ASDH-SDDP in~\cite{Street-Valladao-Lawson-Velloso:2020}).
Contrarily to~\cite{Street-Valladao-Lawson-Velloso:2020} and~\cite{Valladao-Silva-Poggi:2019},
there are no assumption of continuous decision variables, convex costs and constraints, and linear dynamics
up to~\cite[Sect.~4 (Proposed algorithm)]{dowson2020policy}. In that, we share some proximity with~\cite{dowson2020policy}.
For the rest, \cite[Sect.~4 (Proposed algorithm)]{dowson2020policy} deals with the HD information structure for problems that resort to linear programming.

All in all, we have little in common with~\cite{dowson2020policy}.
The author indeed evokes DH and HD --- but the setting is not as formal as ours, and not as general ---
and does not consider DHD, which is our main object of study.
We provide Bellman equations --- whose derivation is not given in the paper, as
it is lengthy, but that can be established from~\cite{Carpentier-Chancelier-DeLara-Martin-Rigaut:2023} ---
for general decision variables (continuous, discrete or mixed),
costs and constraints (convex or not), and dynamics (linear or not).
In~\cite[\S3.4]{dowson2020policy}, a HD Bellman equation~\cite[Equation~(3)]{dowson2020policy}
and a DH Bellman equation~\cite[Equation~(4)]{dowson2020policy} are stated.
Finally, we tackle the issue of DHD with two timescales, whereas
  \cite{dowson2020policy} deals with a single timescale.

 \subsubsection{Synthesis of our contributions in comparison with the literature}
 \label{Synthesis_of_our_contribution_in_comparison_with_the_literature}

 In Table~\ref{tab:Synthesis_of_our_contribution_in_comparison_with_the_literature},
 we highlight the main features of our contributions in comparison with \cite{Street-Valladao-Lawson-Velloso:2020,Valladao-Silva-Poggi:2019}.

\begin{table}[!ht]
  \centering
  \begin{tabular}{|c||c|c|c|c|}
    \hline
    features
    & {\tiny \cite{Street-Valladao-Lawson-Velloso:2020}}
    & {\tiny \cite{Street-Valladao-Lawson-Velloso:2020}}
    & {\tiny \cite{Valladao-Silva-Poggi:2019}}   & our paper
    \\
    &  {\tiny fullDH} &  {\tiny  ASDH} &  &  \\\hline \hline
    general costs, constraints, dynamics  &&&& \checkmark
    \\ \hline
    binary decision variables &&&& \checkmark
    \\ \hline
    SDP &&&& \checkmark
    \\ \hline
    SDDP &\checkmark &\checkmark &\checkmark &
    \\ \hline
    no state extension &\checkmark &&& \checkmark
    \\ \hline
    two timescales &&&& \checkmark
    \\ \hline
  \end{tabular}
  \caption{Synthesis of our contributions in comparison with the literature
    \label{tab:Synthesis_of_our_contribution_in_comparison_with_the_literature}}
\end{table}

\subsection{Structure of the paper}
\label{Structure_of_the_paper}

The paper is organized as follows.

In Sect.~\ref{section:AdequacyProblem}, we present the physical modelling of the problem. First we introduce the
 two-timescale timeline, made of hours (because of hourly energy balance constraints) and weeks (because of weekly
  planning of decisions).
 Second, we define the physical variables needed to describe the system operation
  as well as the linking constraints and the
  economical cost functions.

In Sect.~\ref{section:AdequacyProblemHD}, once the physical model is set, we move on to the mathematical
 formulation of the adequacy resource problem
as a multistage stochastic optimisation problem, focusing on the current practice at RTE for the
information structure modelling, that is, \emph{weekly hazard-decision}.

  In Sect.~\ref{section:AdequacyProblemDHD}, we present a \emph{decision-hazard-decision}  structure that
  considers both ``here-and-now" and ``wait-and-see" decisions in the context of prospective studies.
  The decision-hazard-decision structure is used to solve the adequacy problem with two timescales. The decision stages in the decision-making
process are separated into a nonanticipative  planning step and a recourse corrective step.
As a consequence, the resolution of each stage in the multistage stochastic optimisation problem becomes a
two-stage problem in which the first stage decisions (\emph{slow} decisions) are made before knowing the uncertainties,
and the  second stage or recourse decisions (\emph{fast} decisions) are made once the weekly block of uncertainties is
 known.
It can be interpreted that the slow decisions are associated with
the unit commitment step and that the fast decisions are associated with the unit modulation.
We obtain as a result a problem formulation that improves the information model by being less
anticipative but still allows us to apply temporal decomposition methods.
Once the new information structure is described, we present the corresponding mathematical formulation of the
problem  and the associated Bellman equations giving the usage values.

In Sect.~\ref{section:numericalstudy}, we present numerical results for a  case study comparing both information
structures: \emph{hazard-decision} and \emph{decision-hazard-decision} in the context of prospective studies.
 We observe that the choice of the information structure when computing usage values can modify the merit order in
the system, that is,  the order of the storage usage values with respect to the thermal units prices.
We aim to illustrate the effect of modelling the information structure in a simple case,
and it is not our goal to discuss possible algorithmic approaches to improve performance when computing usage values.
We address the resolution of the two-stage problem at each stage of the multistage optimisation problem in its
 extensive formulation, knowing that more sophisticated techniques can be applied to scale the study. In this paper
 we focus in one storage (or few storages); to address the issue of scalability, we will count on decomposition methods (rather than SDDP)
 like in \cite{Pacaud-DeLara-Chancelier-Carpentier:2022}.

Finally, in Sect.~\ref{Section:Conclusion},  we conclude on the relevance of the choice of information
structures in the computation of usage values, hence on the optimal allocation of resources.

In Appendix~\ref{Appendix:Hourlycomposition}, we detail the
hourly composition of the storage dynamics. 
In Appendix~\ref{Appendix:HourlyRecourseBellman}, we present the ideal
\emph{decision-hazard-decision information structure with hourly recourse}
discussed in~\S\ref{Our_contribution}.

\section{Physical and economical model of the energy system}
\label{section:AdequacyProblem}
In this Sect.~\ref{section:AdequacyProblem}, we describe the physical and economical model of the energy system we consider.
In \S\ref{subsection:timeline}, we present the timeline with two timescales. In \S\ref{subsection:variablesdefinition}, we
define the variables to model the system. In \S\ref{subsection:constraints}, we introduce the system dynamics and energy balance.
Finally, in \S\ref{subsection:Costfunctions}, we present the cost function modelling.

\subsection{Timeline definition}
\label{subsection:timeline}
We consider a timeline with a long timescale and a short timescale.
The short and long timescales could be any two scales, as long as one is larger than the other.
In this work, the long timescale is given by weeks that are represented by a finite totally ordered set
$\np{\WEEK, \preceq}$,
where $\week\successor$ is the successor of $\week\in \WEEK$  and $\week\predecessor$
its predecessor: \(\binf{\week} \prec \cdots \prec \week\predecessor \prec \week
\prec \week\successor \prec \cdots \prec \bsup{\week}\) (where $\prec$ is the strict order associated with the order $\preceq$).
Then, $\WEEK = \nc{\binf{\week}, \bsup{\week}}$.
The short timescale, hours in this case, is represented by a finite totally ordered set $\np{\HOUR, \preceq}$:
\(\binf{\hour} \prec  \cdots  \prec \hour\predecessor \prec\hour \prec\hour\successor \prec \cdots
\prec \bsup{\hour}\).
Then, $\HOUR = \nc{\binf{\hour}, \bsup{\hour}}$.

 To unify the timescale we consider the product set $\timeline$ ordered as follows:
\begin{align}
  \label{eq:timeline}
  \winfhinf&
          \prec\cdots\prec \np{\week\predecessor,\bsup{\hour}}
          \prec \np{\week, \binf{\hour}}
          \prec \np{\week,\binf{\hour}\successor}\prec \cdots
          \nonumber\\
  \cdots &
           \np{\week, \bsup{\hour}\predecessor}\prec\np{\week,\bsup{\hour}}
           \prec \np{\week\successor,\binf{\hour}}\prec \cdots \prec \wsuphsup
           \eqfinp
\end{align}
We consider a period of one year, and $\winfhinf$ is the instant corresponding
to the first hour of the first week of the year, and
$(\bsup{\week},\bsup{\hour})$ is the last hour of the last week
of the year.
We need to define an extra time $\whlast$ at its
end to handle the resulting state of the last decision.
The extended unified timeline $\extimeline$ is defined as $\timeline \cup
\na{\whlast}$.

We define  \(\closedWopen{\week} =\bp{ {\whinf}, {\whinfsuccessor}, \dots,  {\whsup}} \)
and  \( \allowbreak \openWclosed{\week} = \left({\whinfsuccessor}, \dots,  {\whsup},\right.\) \(\left.{\wsuccessorhinf}\right) \).
Thus, we use a simple bracket $\left[\right.$ or $\left.\right]$ to denote intervals of the
elementary timelines $\np{\HOUR, \preceq}$ and  $\np{\WEEK, \preceq}$. By contrast, we use
double brackets $\llbracket$ or $\rrbracket$  for the composite (product) timeline $\np{\extimeline, \preceq}$.

The different possibilities to index a variable (respectively a function) by time are detailed in
Table~\ref{Table:variablesnotationresume}
(respectively in  Table~\ref{table:functionsnotationresume}).
\begin{table}[!ht]
  \renewcommand{\arraystretch}{1.5}
  \centering
  \begin{tabular}{|c|c|c|}
\hline
                      Index                                    &
                      Notation                            &
                      Description
  \\ \hline \hline
  \multicolumn{1}{|c|}{$\wh$}                                         &
                      $z_{\wh}$                                       &
                      Variable at $\wh$
  \\ \hline
  \multicolumn{1}{|c|}{$\week$}                                       &
                      $z_{\week}$                                     &
          \begin{tabular}[c]{@{}c@{}}
                      Representative variable for the week $s$ \\ corresponding to
                      the variable at $\whinf$
          \end{tabular}
  \\ \hline
  \multicolumn{1}{|c|}{$\closedWopen{\week}$}                         &
                      $z_{\closedWopen{\week}}$                       &
          \begin{tabular}[c]{@{}c@{}}
            Sequence of hourly variables given by \\
                        $\bp{ z_{\whinf}, z_{\whinfsuccessor}, \dots,  z_{\whsup}}$
          \end{tabular}
  \\ \hline
  \multicolumn{1}{|c|}{$\openWclosed{\week}$}                         &
                      $z_{\openWclosed{\week}}$                       &
          \begin{tabular}[c]{@{}c@{}}
            Sequence of hourly variables given by \\
                      $\bp{z_{\whinfsuccessor}, \dots,  z_{\whsup}, z_{\wsuccessorhinf}}$
          \end{tabular}
  \\ \hline
  \end{tabular}
  \caption{Variables notation
    \label{Table:variablesnotationresume}}
\end{table}

\begin{table}[!ht]
  \renewcommand{\arraystretch}{1.5}
    \centering
\begin{tabular}{|c|c|c|}
\hline
                   Index                                     &
                  {Notation}                             &
                  {Description}
\\ \hline \hline
\multicolumn{1}{|c|}{$\wh$}                                         &
                  $\phi_{\wh}$                                       &
                  Function expression at $\wh$
\\ \hline
\multicolumn{1}{|c|}{$\week$}                                       &
                  $\phi_{\week}$                                     &
      \begin{tabular}[c]{@{}c@{}}
                  Characteristic aggregation of the hourly
                  functions  \\ $\phi_{\wh}$ for
                  the week $s$
      \end{tabular}
\\ \hline
\multicolumn{1}{|c|}{
                         $\closedWopen{\week}$
                          } &
           $\phi_{\closedWopen{\week}}$                                   &
      \begin{tabular}[c]{@{}c@{}}
        Sequence of hourly functions $\phi_{\wh}$ given by \\
        $\bp{ \phi_{\whinf}, \phi_{\whinfsuccessor}, \dots,  \phi_{\whsup}}$
      \end{tabular}
\\ \hline
\multicolumn{1}{|c|}{
                         $\openWclosed{\week}$ } &
                         $\phi_{\openWclosed{\week}}$                            &
      \begin{tabular}[c]{@{}c@{}}
        Sequence of hourly functions $\phi_{\wh}$ given by \\
        $\bp{ \phi_{\whinfsuccessor},  \dots, \phi_{\whsup}, \phi_{\wsuccessorhinf}}$
      \end{tabular}
\\ \hline
\end{tabular}
\caption{Functions notation
\label{table:functionsnotationresume}}
\end{table}
The characteristic aggregation in Table~\ref{table:functionsnotationresume} could be a sum over~$\hour\in\HOUR$, a composition
with respect to the state, or a combination of both (see Appendix~\ref{Appendix:Hourlycomposition} for further details).

\subsection{Physical variables}
\label{subsection:variablesdefinition}
The following is a description of the system components.
We classify the variables accordingly to their type: decision (in the hand of the decision-maker),
 uncertainty (exogenous),
state (storage) and slack (energy not supplied in the system).

\subsubsection{Thermal units modelling}
We consider a thermal fleet composed of thermal units whose variables are detailed in
 Table~\ref{Table:thermalvariables}
for the units indexed by $i\in\II$.
\begin{table}[!ht]
  \renewcommand{\arraystretch}{1.5}
  \centering
  \begin{tabular}{|c|c|c|}
  \hline
  Description &  Type & Notation \\
  \hline  \hline
  On/Off & Decision & $\techcommitmenti \in \na{0,1}  $  \\ \hline
  Modulation    & Decision & $   \tech^{i}\in \na{0} \cup \bc{\binf{\tech}^i, \bsup{\tech}^i}  $  \\ \hline
  Availability & Uncertainty & $\techuncertaini  \in \na{0,1}$ \\ \hline
  \end{tabular}
  \caption{Thermal units variables
  \label{Table:thermalvariables}}
\end{table}
The ``On" decision is associated with $\techcommitmenti=1$ and the ``Off" decision
with $\techcommitmenti=0$. Observe that the decision $\techcommitmenti$ is taken at each hour and does not represent a change in the state of the unit.
The decision $ \tech^{i}$ denotes the power modulation once the unit is on.
The availability of the thermal units is modelled with the (uncertainty) variable~$\techuncertaini$:
 when it is equal to~$0$ the unit is not available to use
  and, when it is equal to~$1$, the unit is available to use.

The collections of on/off decisions, modulation decisions and availabilities variables
of all thermal units are denoted by $\techcommitment~=~\nseqp{\techcommitmenti}{i\in\II}$,
$\tech~=~\nseqp{\tech^i}{i\in\II}$ and $\techuncertain~=~\nseqp{\techuncertaini}{i\in\II}$ respectively.

\subsubsection{Storage modelling}

In Table~\ref{Table:storagevariables}, we introduce the variables related to the storage management.
We consider different variables for pumping $(\pumping)$ and turbining $(\turbining)$ decisions
so that we take into account the
pumping efficiency in the storage. The variable~$\stock$ denotes the level in the storage (stock), that
is, the physical state of the storage.
\begin{table}[!ht]
  \renewcommand{\arraystretch}{1.5}
  \centering
  \begin{tabular}{|c|c|c|}
  \hline
  Description &  Type & Notation \\
  \hline  \hline
  Pumping & Decision & $\pumping \in \bc{0, \bsuppumping}  $  \\ \hline
  Turbining    & Decision & $\turbining \in \bc{0, \bsupturbining}  $  \\ \hline
  Storage level & Physical state & $\stock \in \bc{\binf{\stock}, \bsup{\stock}} $\\ \hline
  \end{tabular}
  \caption{Storage variables
  \label{Table:storagevariables}}
\end{table}

\subsubsection{Residual demand modelling}
The \emph{residual demand} ($\uncertain^{\demand} $) is the difference
between demand and non-dispatchable production. This allows to group in one variable several uncertainty sources
 such as  the wind production, the solar production, the demand, etc. The variable $\uncertain^{\demand} $
 is classified as uncertainty.

\subsubsection{Slack variables}
We introduce a variable $\ens$, classed as \emph{slack}, to model the \emph{energy not supplied} in the system.
This variable will appear in the forthcoming energy balance~\eqref{eq:BalanceEquations} and
 cost function~\eqref{eq:CostFunctions}.

\subsection{System dynamics and energy balance}
\label{subsection:constraints}

We present now the linking constraints between variables.

 \subsubsection{Storage dynamics}
The dynamics function
\label{eq:StorageDynamics}
\begin{align}
  \label{eq:StorageHourlyDynmics}
    \dynamics\np{ \stock, \pumping ,\turbining
           }
 = \stock + \eta \pumping - \turbining
\end{align}
 describes the evolution of the storage level --- as a function of the current storage level~$\stock$, 
of the pumping~$\pumping$ and turbining~$\turbining$ decisions --- from one short
 time stage to the next. The parameter $\eta\in\nc{0,1}$ is the pumping efficiency
 of the storage.

The hourly dynamics $\dynamics$ induces a weekly temporal coupling linking the storage level at the beginning of a week
$\week$ with the level at the beginning of the following week $\week\successor$.
Therefore, we also consider a weekly dynamics $\dynamics_{\week}$
  given by the hourly composition of the dynamics in~\eqref{eq:StorageHourlyDynmics}.
  The composition is detailed in Appendix~\ref{Appendix:Hourlycomposition}.
   It gives,
as a result, the storage level $\stock_{\wsuccessorhinf}$ at the beginning of the following week,
by summing the total difference between pumping (positive taking into account its efficiency) and turbining (negative) during the week
  to the storage level at the beginning of the current week $\stock_{\whinf}$.



\subsubsection{Thermal units' production output}

The effective output of the unit is constrained by its availability and the on/off decision.
More precisely, the output production $\techproduction^i$ of the $i$-th thermal unit
  not only depends on the decision~$\tech^{i}$, but also on the on/off decision
$\techcommitmenti$ and on the availability uncertainty $\techuncertaini$ as follows
\begin{align}
\techproduction^i \np{\techcommitmenti,
                          \tech^i,
                           \techuncertaini} =
     \tech^i \times  \min\na{\techcommitmenti, \techuncertaini}
    \eqfinp
\end{align}
As defined, the production of the $i$-th unit belongs to the same set as $\tech^i$, that is,
\begin{align}
 \techproduction^i\np{\techcommitmenti,
                          \tech^i,
                           \techuncertaini} \in \na{0} \cup \bc{\binf{\tech}^i, \bsup{\tech}^i}  \eqfinv
\end{align}
taking the value zero whenever the off
decision is made ($\techcommitmenti= 0$) or the unit is not available ($\techuncertaini=0$).
The collection of output production of all thermal units is
 denoted by $\techproduction~=~\nseqp{\techproduction^i}{i\in\II}$.

\subsubsection{Balance equation}

Satisfying the energy balance equation is the main goal of the adequacy resource problem. That is,
at every hour, the total energy production in the system should be equal to the total energy consumed in the system  (which includes the pumping).
Ideally, the energy balance constraint is written as the equality
\begin{align}
  \label{eq:equalityBalance}
  \overbrace{\turbining +
  \sum_{i\in\II} \techproduction^i }^{\substack{\text{total}\\\text{production}}} +
  \overbrace{\ens}^{\substack{\text{energy not}\\\text{supplied}}}
    =
    \overbrace{  \pumping +
  \uncertain^{\demand} }^{\substack{\text{pumping demand +}\\\text{residual demand}}}
  \eqfinp
\end{align}

We rather formulate this balance equation as an inequality constraint to avoid infeasibility problems
 due to the thermal units' minimum power constraints: indeed, infeasibility could happen when less energy is required than
 the minimum power of the last (or more expensive) unit on to meet the demand.
If the cause of the balance infeasibility is the lack of available production, the slack variable $\ens$ will
take positive values, measuring how far from meeting the demand the system is. Thus, the hourly balance equation is
given by
\begin{subequations}
  \label{eq:BalanceEquations}
\begin{align}%
  \label{eq:hourlyBalance}
  \nodeBalance&(\turbining,\pumping,
            \techproduction ,
                  \uncertain^{\demand},
                    \ens)
                    \geq 0 \eqfinv
                    \\
  \label{eq:hourlyBalanceExplicit}
\text{with } \nodeBalance \bp{\turbining,\pumping,
  \techproduction,
  \uncertain^{\demand},
    \ens} &=
\Bp{\turbining +
  \sum_{i\in\II} \techproduction^i +
  \ens }
  - \bp{
  \pumping +
  \uncertain^{\demand}}
  \eqfinp
\end{align}
\end{subequations}


\subsection{Cost functions}
\label{subsection:Costfunctions}

We now introduce the hourly cost
function~$\instantaneouscost\np{\techcommitment, z^{\mathsf{th}},  \techproduction,\techuncertain, \ens} $
as
\begin{align}
\label{eq:CostFunctions}
  \instantaneouscost
    \np{
      \techcommitment,
      z^{\mathsf{th}},
      \techproduction,
      \techuncertain,
        \ens}
  =
  \sum_{i\in\II} \Bigl(
      \startupcost^i &\times \max\ba{\techcommitmenti- z^{\mathsf{th},i}, 0}
  +
    \variablecost^i \times \techproduction^i
    \Bigr)  \nonumber
  +
    \penalizationCost \times\ens
  \eqfinp
\end{align}

The cost of meeting the demand is the operating cost of the thermal units.
We model the thermal cost for each unit within the sum over $i\in\II$ in Equation~\eqref{eq:CostFunctions} with two components. 
The first one is associated with decision of
switching on a unit, that is, when at one hour $\techcommitmenti =0$
and in the following hour  $\techcommitmenti =1$. The second component corresponds to the variable cost, and is
proportional to the power modulation~$\tech^i$ of the unit.
In addition, we model the penalization on the energy not supplied as a cost.
This penalization cost is much higher than the thermal units' cost to ensure that the energy demand is not
provided only in cases where there is no other solution.
 
The parameters $\startupcost^i$ and $\variablecost^i$ correspond to the unit's start-up
and modulation costs, whereas $\penalizationCost$ is the penalization parameter for the not supplied energy in the system.
The variable $z^{\mathsf{th}}$ is introduced to take into account a temporal shift of
the on/off decision $\techcommitmenti$ since the start up cost
is associated with a change from off to on between two consecutive hours
(see Tables~\ref{Table:HDmathematicalnotation} and \ref{Table:DHDmathematicalnotation}).


The weekly cost $\instantaneouscost_{\week}$ will be defined as
the sum of the hourly costs $\instantaneouscost$ within the week in Tables~\ref{Table:HDmathematicalnotation} and
\ref{Table:DHDmathematicalnotation}.


Up to now, we have introduced the physical and economical modelling of the problem.
In Sect.~\ref{section:AdequacyProblemHD},
we will present the  current practice for the adequacy problem mathematical
formulation, to compute usage values, focusing on the  information structure modelling.

\section{Current information modelling in weekly hazard-decision}
\label{section:AdequacyProblemHD}

When modelling a stochastic multistage optimization problem, it is necessary to define the
information structure, that is, a model that describes the information available at
each stage of the decision-making process.

In \S\ref{subsection:NotationRandomVariables}, we introduce notations for random variables
and measurability constraints. In \S\ref{subsection:WeeklyHD}, we present the weekly hazard-decision information structure.
In \S\ref{subsection:MultistageStochasticOptHD}, we formulate the multistage stochastic optimization problem.
In \S\ref{subsection:BellmanHD}, we deduce the corresponding Bellman equations in hazard-decision.

\subsection{Notation for random variables and measurability constraints}
\label{subsection:NotationRandomVariables}

To formulate the multistage stochastic optimization problem,
we model the uncertainties as random variables and, as a consequence, the states and
controls are random variables as well.
For this purpose, we consider a probability space~$(\Omega,\trib,\prbt)$.
Random variables are measurable functions from the measurable space~$(\Omega,\trib)$
towards some~$\RR^n$, equipped with its Borel \sigmafield~$\borel{\RR^n}$,
and will be denoted by bold letters like~$\va{Z}$.
Then, the information structures are mathematically modelled as measurability constraints.
We say that $\va{Z_1}$ is measurable with respect to $\va{Z_2}$
if~$\nsigmaf{\va{Z_1}} \subset \nsigmaf{\va{Z_2}}$, that is, the \sigmafield generated
by~$\va{Z_1}$ is included in (less rich than) the \sigmafield generated by~$\va{Z_2}$. 
Thanks to Doob Theorem {\cite{Dellacherie-Meyer:1975}}, measurability constraints are equivalently expressed
by means of functions as follows
\begin{align}
    \nsigmaf{\va{Z_1}}\subset \nsigmaf{\va{Z_2}}   \iff
  \text{there exists a measurable function } \varphi
                     \text{ such that  } \va{Z_1} =
                        \varphi\np{\va{Z_2}} 
                        \eqfinp
\end{align}
Practically, the function $\varphi$ is what we call policy or strategy,
when $\va{Z_2}$ represents the information disclosed when making a decision.

\subsection{Weekly hazard-decision information structure (HD)}
\label{subsection:WeeklyHD}
In this Sect.~\ref{section:AdequacyProblemHD}, we consider information structures that assume a \emph{weekly} disclosure
of the information: once the uncertainty of the first hour of the week is known,
 the whole collection of uncertainties for the week is also known.
The weekly vectors presented hereafter follow the notation
 given in Table~\ref{Table:variablesnotationresume}.
 Let \( \uncertain_{\openWclosed{\week}} =
 \np{\uncertain^{\demand}_{\openWclosed{\week}},\techuncertain_{\openWclosed{\week}}}\)
  be the vector composed of the weekly demand and weekly thermal units availabilities.
In the same way, the vector composed of the collection of weekly controls is denoted by
$\recoursecontrol_{\openWclosed{\week}} =\np{\techcommitment_{\openWclosed{\week}},  \tech_{\openWclosed{\week}},
\pumping_{\openWclosed{\week}}, \turbining_{\openWclosed{\week}},
\ens_{\openWclosed{\week}}}$.

The current practice to model the information structure
is a \emph{hazard-decision} $(\HD)$ structure in the weekly timescale as illustrated
in Fig.~\ref{Figure:WeeklyHD}.

\begin{figure}[htbp]
  \centering
  \resizebox{0.8\columnwidth}{!}{
       \begin{tikzpicture}
        \filldraw[black] (-7.,0) circle (2pt) node[anchor=north] {$\winfhinf$};
        \draw[dashed, thick] (-7.,0) -- (-5.25,0);
        \filldraw[black] (-5.25,0) circle (2pt) node[anchor=north] {$\whinf$};
        \draw[black, thick] (-5.25,0) -- (-4,0);
        \filldraw[black] (-4,0) circle (2pt) node[anchor=north] {$\whinfsuccessor$};
        \draw[dashed, thick] (-4,0) -- (-2.25,0);
        \filldraw[black] (-2.25,0) circle (2pt) node[anchor=north] {$\wh\predecessor$};
        \draw[black, thick] (-2.25,0) -- (-1,0);
        \filldraw[black] (-1,0) circle (2pt) node[anchor=north] {$\wh$};
        \draw[black, thick] (-1,0) -- (0.25,0);
        \filldraw[black] (0.25,0) circle (2pt) node[anchor=north] {$\wh\successor$};
        \draw[dashed, thick] (0.25,0) -- (2,0);
        \filldraw[black] (2,0) circle (2pt) node[anchor=north] {$\whsup$};
        \draw[dashed, thick] (0.25,0) -- (2,0);
        \draw[black, thick] (02 ,0) -- (3.25,0);
        \filldraw[black] (3.25,0) circle (2pt) node[anchor=north] {$\wsuccessorhinf$};
        \draw[-latex] (3.25,1)-- (3.25,0.15);
        \node (vs) at (3.25,1.3) {\footnotesize {${\recoursecontrol_{\openWclosed{\week}}}$}};
        \draw[dashed, thick] (3.25,0) -- (4.5,0);
    \draw [decorate,
    decoration = {brace, raise=5pt, mirror,
        amplitude=8pt}] (-4,-0.5) --  (3.25,-0.5);
    \node (uncertain)[align=center] at (-.375, -1.2) {\footnotesize	 Uncertainty $\uncertain_{\openWclosed{\week}}$ of the week $\week$};
    \node () at (-.375, -1.5){\footnotesize	 $\bp{\uncertain_{\whinfsuccessor}, \dots, \uncertain_{\wh}, \dots, \uncertain_{\wsuccessorhinf}}$};
    \draw[ -latex,  dashed] (3.25 ,1)-- (-4, 0.15);
    \draw[ -latex,  dashed] (3.25 ,1)-- (-2.25,0.15);
    \draw[ -latex,  dashed] (3.25 ,1)-- (-1,0.15);
    \draw[ -latex,  dashed] (3.25 ,1)-- (0.25,0.15);
    \draw[ -latex,  dashed] (3.25 ,1)-- (2,0.15);
    \draw[ -latex,  dashed] (3.25 ,1)-- (3.25,0.15);
    \draw [decorate,
           decoration = {brace, raise=5pt,
            amplitude=8pt}] (-7,1.4) --  (3.25,1.4);
        \node (W-HD) at (-2.125,2.1) { \footnotesize	Span of available information before making decision $\recoursecontrol_{\openWclosed{\week}}$};
        \draw[dashed, thick] (3,0.8) -- (3,0.8);
                \draw[dashed, thick] (-7.,1.4) -- (-7.,0);
\end{tikzpicture}
  }
  \caption{Weekly hazard-decision information structure.
      An arrow maps the available information towards the decision,
  so that, here, the information structure is partly anticipative as some arrows go from the right
  to the left.
\label{Figure:WeeklyHD}}
\end{figure}
The collection \( \recoursecontrol_{\openWclosed{\week}} =
 \bp{\recoursecontrol_{\whinfsuccessor}, \dots, \recoursecontrol_{\wsuccessorhinf}}\)
 of hourly decisions for the week  is made
once the block \(\uncertain_{\openWclosed{\week}} =
\bp{\uncertain_{\whinfsuccessor}, \dots, \uncertain_{\wsuccessorhinf}}\)
of uncertainties for the week  is disclosed. In other words,
 when making the decisions for any hour of the week, the demand and the availability of thermal units for
 every hour in the week are already known (in advance). We can also interpret this structure as if the hourly decisions can wait
 until knowing all the uncertainties of the week to be made. In this context, all the hourly decisions within the week are
\emph{one week ahead anticipative} since, when making them, the uncertainties until the end of the week are already known.

Finally, for the weekly hazard-decision structure, the information constraint is represented by the following measurability constraints
\begin{align}
  \sigma\np{\va{\RecourseControl}_{\openWclosed{\week}} } \subset
  \sigma\np{\va{\Uncertain}_{\winfhinf},
              \va{\Uncertain}_{\openWclosed{\binf{\week}}}, \dots,
              \va{\Uncertain}_{\openWclosed{\week}} }
              \eqsepv
              \forall \week \in \WEEK
              \eqfinp
\end{align}

\subsection{Multistage stochastic optimization problem formulation}
\label{subsection:MultistageStochasticOptHD}

In Sect.~\ref{section:AdequacyProblem}, we have presented the physical and economical model of the energy system.
With the new notation in \S\ref{subsection:WeeklyHD}, we present in Table \ref{Table:HDmathematicalnotation}
the corresponding expressions and their compact mathematical versions.
\begin{table}[!h]
  \renewcommand{\arraystretch}{1.2}
  \centering
  \begin{tabular}{c|@{\hspace{0cm}}c@{\hspace{0.1cm}}|@{\hspace{0cm}}c@{\hspace{0.1cm}}|}
   \cline{2-3}
  &  \multicolumn{2}{c|}{Notation}\\
    \hline
    \multicolumn{1}{|c|}{              Function }                                                                                       &
                          Mathematical  &
                      Physical-Economical                                                                                                                     \\
                                                                                                                  \hline
                                                                                                                  \hline
  \multicolumn{1}{|c|}{
                        \begin{tabular}[c]{@{}c@{}}
                          \scriptsize        Storage \\
                          \scriptsize         dynamics
                    \end{tabular}
                            }                                         &
                            \hspace{0.01cm} $\scriptstyle \Dynamics_{{\week}}\bp{\stock_{\whinf}, \recoursecontrol_{\openWclosed{\week}}}$         &
                     $ \scriptstyle \dynamics_{\week}   \bp{ \stock_{\whinf}, \pumping_{\openWclosed{\week}} ,\turbining_{\openWclosed{\week}}}$
                                                                                                                             \\ \hline
  \multicolumn{1}{|c|}{
                    \begin{tabular}[c]{@{}c@{}}
                          \scriptsize       Thermal \\
                          \scriptsize        production
                    \end{tabular}
                            }                                         &
                            \hspace{0.01cm}  $\scriptstyle \Techproduction_{\openWclosed{\week}}\bp{\uncertain_{\openWclosed{\week}}, \recoursecontrol_{\openWclosed{\week}}}$                       &
             $\scriptstyle \Bigl\{\techproduction_{\wh\successor} \bp{\techcommitment_{\wh\successor}, \tech_{\wh\successor},\techuncertain_{\wh\successor}}\Bigr\}_{\hour\in\HOUR} $
  \\ \hline
  \multicolumn{1}{|c|}{
                        \begin{tabular}[c]{@{}c@{}}
                  \scriptsize        Energy \\
                  \scriptsize       balance
                    \end{tabular}
                            }                                         &
                \hspace{0.01cm}   $\scriptstyle \NodeBalance_{\closedWopen{\week}}\bp{{\uncertain_{\openWclosed{\week}}, \recoursecontrol_{\openWclosed{\week}}}} $
                    &
          \begin{tabular}[l]{@{\hspace{0cm}}p{4.5cm}@{\hspace{0cm}}}
              \raggedright $\scriptstyle  \biggl\{\nodeBalance \Bigl( \turbining_{\wh\successor},\pumping_{\wh\successor},  $
                  \\
              \centering   $\scriptstyle  \techproduction_{\wh\successor} \bp{\techcommitment_{\wh\successor},\tech_{\wh\successor}, \techuncertain_{\wh\successor}},$
                 \\
               \raggedleft $\scriptstyle   \uncertain^{\demand}_{\wh\successor},\ens_{\wh\successor} \Bigr)\biggr\}_{\hour\in\HOUR}$
          \end{tabular}
  \\ \hline
  \multicolumn{1}{|c|}{
                        \begin{tabular}[c]{@{}c@{}}
                          \scriptsize          Weekly \\
                          \scriptsize       cost
                    \end{tabular}
                            }                                         &
                  \hspace{0.01cm} $\scriptstyle \InstantaneousCost_{\week}\bp{{\uncertain_{\openWclosed{\week}}, \recoursecontrol_{\openWclosed{\week}}}}$
                     &
            \begin{tabular}[l]{@{\hspace{0cm}}p{4.5cm}@{\hspace{0cm}}}
              \raggedright  $   \scriptstyle \instantaneouscost_{\week} \bp{
               \techcommitment_{\openWclosed{\week}},
           \techproduction_{\openWclosed{\week}}, \techuncertain_{\openWclosed{\week}}, \ens_{\openWclosed{\week}}} = $
           \\
               \raggedleft $
                 \scriptstyle\instantaneouscost \bigl(\techcommitment_{\whinf\successor},
                                                                                0,
                                                                                \techproduction_{\whinf\successor},
                                                                                \techuncertain_{\whinf\successor},
                                                                                 \ens_{\whinf\successor}
                \bigr)+
           \displaystyle \sum_{\scriptscriptstyle \hour\in\HOUR \setminus\na{\binf{\hour}}}
           \scriptstyle\instantaneouscost \bigl(\techcommitment_{\wh\successor},
                                                                      \techcommitment_{\wh},
                                                                      \techproduction_{\wh\successor},$
             \\
              \raggedleft                      $\scriptstyle \techuncertain_{\wh\successor},   \ens_{\wh\successor} \bigr)$
          \end{tabular}
  \\ \hline
  \end{tabular}
  \caption{Correspondence between mathematical and physical and economical notations (Sect.~\ref{section:AdequacyProblem})
    in the weekly hazard-decision framework
  \label{Table:HDmathematicalnotation}}
\end{table}

The weekly cost function $\instantaneouscost_{\week}$ is obtained as the sum of the hourly cost functions within the week.
When computing the hourly cost at $\wh$ with Equation~\eqref{eq:CostFunctions}, the variable $z^{\mathsf{th}}$
correspond to the on/off decision $\techcommitment_{\wh\predecessor}$ at the previous hour.
For simplicity, we neglect the temporal coupling of thermal units from the last hour of one week to the
first hour of the next week, fixing the initial operating state of the units at the beginning of each week to off.
  This modelling choice is made to avoid an augmentation of the dynamic programming state's dimension for each
 thermal unit.
This explains the presence of a zero as second entry in the first term
in the expression of the weekly cost in the last line of Table~\ref{Table:HDmathematicalnotation}.

Considering the definitions in Table~\ref{Table:HDmathematicalnotation}, we formulate the physical
adequacy problem as a stochastic multistage optimization  problem using the weekly hazard-decision information structure:
\begin{subequations}
  \label{eq:WeeklyHDProblemFormulation}
  \begin{align}
    &
      \min_{\va{\Stock}, \va{\RecourseControl}}
        \EE  \biggl[
              \sum_{\week\in \WEEK}
              \InstantaneousCost_{\week}\np{\va{\Stock}_{\whinf},
                          \va{\RecourseControl}_{\openWclosed{\week}},
                          \va{\Uncertain}_{\openWclosed{\week}}}
                  +
                  \FinalCost\np{\va{\Stock }_{\whlast}}			
                  \biggr]
                  \\
      & \suchthat \eqsepv
      \forall \week \in \WEEK
      \nonumber
      \\
      &
      \va{\Stock}_{\winfhinf} = \va{\Uncertain}_{\winfhinf} \eqfinv
      \\
      &
      \va{\Stock}_{\wsuccessorhinf} = \Dynamics_{\week}\np{\va{\Stock}_{\whinf},
                                                              \va{\RecourseControl}_{\openWclosed{\week}}
                                                           }		
              \eqfinv
      \\
      &
      \NodeBalance_{\closedWopen{\week}}\np{\va{\RecourseControl}_{\openWclosed{\week}},
      \va{\Uncertain}_{\openWclosed{\week}}} \geq     \va{0} \eqfinv \label{eq:balanceHD}\\
      &
      \sigma\np{ \va{\RecourseControl}_{\openWclosed{\week}}}
                          \subset \bsigmaf{
                                      \va{\Uncertain}_{\winfhinf},
                                      \va{\Uncertain}_{\openWclosed{\binf{\week}}},
                                      \dots,
                                      \va{\Uncertain}_{\openWclosed{\week\predecessor}},
                                      \va{\Uncertain}_{\openWclosed{\week}}
                                      }
                    \eqfinp
   \end{align}
  \end{subequations}
The final cost $\FinalCost\np{\va{\Stock }_{\whlast}}$ is used to give value to the energy in the storage at
the end of the yearly period.

\subsection{Bellman equations in hazard-decision}
\label{subsection:BellmanHD}

Defining the weekly state $\state_{\week} = \stock_{\whinf}$ (storage level
  at the beginning of the week), we write the weekly Bellman
 equations~\cite{Bellman:1957}
\begin{subequations}
  \label{eq:HDweeklyBellman}
  \begin{align}
    \label{eq:HDweeklyBellmanA}
    & \nBellman{\weeklast\successor}{\HD}{\state_{\weeklast\successor}} =   \FinalCost\np{\state_{\weeklast\successor}}\eqfinv
    \\
    \label{eq:HDweeklyBellmanB}
    & \nBellman{\week}{\HD}{\statew} =
    {\EE} \bigg[
        \min_{\recoursecontrol_{\openWclosed{\week}}  }
        \InstantaneousCost_{\week}\bp{
            \statew,
           \recoursecontrol_{\openWclosed{\week}},
            \va{\Uncertain}_{\openWclosed{\week}}}
        +
    \bBellman{\week\successor}{\HD}{\Dynamics_{\week}\np{
        \statew,\recoursecontrol_{\openWclosed{\week}}}}
                \bigg]
                \eqfinv
\end{align}
\end{subequations}
where the minimum inside Equation~\eqref{eq:HDweeklyBellmanB} is computed subject to the constraint~\eqref{eq:balanceHD}.

If the sequence $\bp{\va{\Uncertain}_{\openWclosed{\binf{\week}}}, \dots,
\va{\Uncertain}_{\openWclosed{\week}} , \dots, \va{\Uncertain}_{\openWclosed{\bsup{\week}}}}$
of uncertainties is weekly  independent,
the weekly Bellman equations provide an optimal solution for Problem~\eqref{eq:WeeklyHDProblemFormulation}.
We highlight that, to get the optimal solution, the hourly uncertainties
  \(\va{\Uncertain}_{\openWclosed{\week}}=  \bp{\va{\Uncertain}_{\whinfsuccessor}, \dots, \va{\Uncertain}_{\wh\successor}, \dots, \va{\Uncertain}_{\wsuccessorhinf}} \)
within the week do not need to be assumed to be independent (from one hour to another).

It is well known that, for all $\week\in \WEEK$, the function $\nBellmanf{\week}{\HD}$ satisfies
\begin{subequations}
  \begin{align}
  \nBellman{\week}{\HD}{\statew} =
  &  \min_{\va{\RecourseControl}_{\openWclosed{\week}},\dots, \va{\RecourseControl}_{\openWclosed{\bsup{\week}}}}
        \EE
        \biggl[
            \sum_{\week'=\week}^{\bsup{\week}}
            \InstantaneousCost_{\week'}\np{\va{\Stock}_{\wprimehinf},
                        \va{\RecourseControl}_{\openWclosed{\week'}},
                        \va{\Uncertain}_{\openWclosed{\week'}}}
                +
                \FinalCost\np{\va{\Stock }_{\whlast}}			
                \biggr]
                \\
    & \suchthat \eqsepv
    \forall \week' \in \nc{\week,\dots,\bsup{\week}}
    \nonumber
    \\
    &
    \va{\Stock}_{\wprimehinf} =  \statew \eqfinv
    \\
    &
    \va{\Stock}_{\wprimesuccessorhinf} = \Dynamics_{\week'}\np{\va{\Stock}_{\wprimehinf},
                                                            \va{\RecourseControl}_{\openWclosed{\week'}}
                                                         }		
            \eqfinv
    \\
    &
    \NodeBalance_{\closedWopen{\week'}}\np{\va{\RecourseControl}_{\openWclosed{\week'}},
    \va{\Uncertain}_{\openWclosed{\week'}}} \geq     \va{0} \eqfinv \\
    &
    \sigma\np{ \va{\RecourseControl}_{\openWclosed{\week'}}}
                        \subset \bsigmaf{
                                    \va{\Uncertain}_{\winfhinf},
                                    \va{\Uncertain}_{\openWclosed{\binf{\week}}},
                                    \dots,
                                    \va{\Uncertain}_{\openWclosed{\wprimepredecessor}},
                                    \va{\Uncertain}_{\openWclosed{\week'}}
                                    }
                  \eqfinv
 \end{align}
\end{subequations}
so that the  value $ \nBellman{\week}{\HD}{\statew}$ of the Bellman function is interpreted as the future optimal cost
when, at week $\week$, the storage level is $\statew$ and the weekly hazard-decision
information structure is considered.

\bigskip
The need for another information structure arises out of the fact that the current
approach is fully anticipative in the week, as illustrated in the minimum inside the expectation
on the right hand side of the Bellman equation~\eqref{eq:HDweeklyBellmanB}.
When assuming that all the uncertainties for the week are known at the moment of making
a decision, we implicitly suppose  that all decisions  are flexible and
can wait until knowing the uncertainties to be made.
But it is known that certain on/off decisions cannot be made instantaneously
and need to be planned in advance.
For this purpose, we introduce another information structure in Sect.~\ref{section:AdequacyProblemDHD}.

\section{Information modelling in weekly decision-hazard-decision}
\label{section:AdequacyProblemDHD}

   In this Sect.~\ref{section:AdequacyProblemDHD}, we present the \emph{weekly decision-hazard-decision} $(\DHD)$
   information structure. As discussed in~\S\ref{Our_contribution},
   this structure is a compromise between the current weekly hazard-decision
   described in Sect.~\ref{section:AdequacyProblemHD}
   and the ideal information structure with weekly planning decisions and hourly recourses
   (detailed in Appendix~\ref{Appendix:HourlyRecourseBellman}).
In the weekly decision-hazard-decision information structure,
 there are decisions that cannot be modelled as (weekly) anticipative.
 We classify the decisions in the system modelling between:
 \begin{itemize}
   \item \emph{planning} or \emph{decision-hazard} (here-and-now) decisions: denoted by~$\controltilde$,
   \item \emph{recourse} or \emph{hazard-decision} (wait-and-see) decisions: denoted by~$\recoursecontroltilde$.
   \end{itemize}
In \S\ref{subsection:WeeklyDHD}, we present and detail the weekly decision-hazard-decision information structure.
In \S\ref{subsection:MultistageStochasticOptDHD}, we formulate the multistage stochastic optimization problem.
In \S\ref{subsection:BellmanDHD}, we give the corresponding Bellman equations in decision-hazard-decision.
In \S\ref{subsection:TheoreticalBellmanComparison}, we compare the Bellman value functions according to
the underlying information structures, namely~$\HD$ and~$\DHD$.

 \subsection{Weekly decision-hazard-decision information structure (DHD)}
\label{subsection:WeeklyDHD}

As shown in Fig.~\ref{Figure:WeeklyDHD},
the planning decisions $\controltilde_{\closedWopen{\week}}$ are made before knowing the uncertainties
$\uncertain_{\openWclosed{\week}}$ for the week (knowing
only the past uncertainties);
then, the weekly block of uncertainties is disclosed, and the corrective actions are made,
that is, the recourse controls $\recoursecontroltilde_{\openWclosed{\week}}$.
\begin{figure}[htbp]
  \centering
  \resizebox{0.8\columnwidth}{!}{
    \begin{tikzpicture}
     \filldraw[black] (-7.,0) circle (2pt) node[anchor=north] {$\winfhinf$};
     \draw[dashed, thick] (-7.,0) -- (-5.25,0);
     \node (us) at (-5.25,1.2) {\footnotesize {${\controltilde_{\closedWopen{\week}}}$}};
     \draw [decorate,
           decoration = {brace, raise=5pt,
               amplitude=8pt}] (-7.,1.25) --  (-5.25,1.25);
           \node (W-DH) [align=left] at (-5.6,2.15) { \scriptsize Information  before \\ \scriptsize making decision $\controltilde_{\closedWopen{\week}}$};
     \filldraw[black] (-5.25,0) circle (2pt) node[anchor=north] {$\whinf$};
     \draw[black, thick] (-5.25,0) -- (-4,0);
     \filldraw[black] (-4,0) circle (2pt) node[anchor=north] {$\whinfsuccessor$};
     \draw[dashed, thick] (-4,0) -- (-2.25,0);
     \filldraw[black] (-2.25,0) circle (2pt) node[anchor=north] {$\wh\predecessor$};
     \draw[black, thick] (-2.25,0) -- (-1,0);
     \filldraw[black] (-1,0) circle (2pt) node[anchor=north] {$\wh$};
     \draw[black, thick] (-1,0) -- (0.25,0);
     \filldraw[black] (0.25,0) circle (2pt) node[anchor=north] {$\wh\successor$};
     \draw[dashed, thick] (0.25,0) -- (2,0);
     \filldraw[black] (2,0) circle (2pt) node[anchor=north] {$\whsup$};
     \draw[dashed, thick] (0.25,0) -- (2,0);
     \draw[black, thick] (02 ,0) -- (3.25,0);
     \filldraw[black] (3.25,0) circle (2pt) node[anchor=north] {$\wsuccessorhinf$};
     \node (vs) at (3.25,1.2) {\footnotesize {${\recoursecontroltilde_{\openWclosed{\week}}}$}};
     \draw[dashed, thick] (3.25,0) -- (4.5,0);
 \draw [decorate,
 decoration = {brace, raise=5pt, mirror,
     amplitude=8pt}] (-4,-0.5) --  (3.25,-0.5);
 \node (uncertain)[align=center] at (-.375, -1.2) {\footnotesize	 Uncertainty $\uncertain_{\openWclosed{\week}}$ of the week $\week$};
 \node () at (-.375, -1.5){\footnotesize	 $\bp{\uncertain_{\whinfsuccessor}, \dots, \uncertain_{\wh}, \dots, \uncertain_{\wsuccessorhinf}}$};
          \draw[ -latex,  dashed] (3.25 ,1)-- (-4, 0.15);
          \draw[ -latex,  dashed] (3.25 ,1)-- (-2.25,0.15);
          \draw[ -latex,  dashed] (3.25 ,1)-- (-1,0.15);
          \draw[ -latex,  dashed] (3.25 ,1)-- (0.25,0.15);
          \draw[ -latex,  dashed] (3.25 ,1)-- (2,0.15);
          \draw[ -latex,  dashed] (3.25 ,1)-- (3.25,0.15);
 \draw [decorate,
        decoration = {brace, raise=5pt,
         amplitude=8pt}] (-7,2.5) --  (3.25,2.5);
     \node (W-HD) at (-1.8,3.1) { \scriptsize	Information before
     making decision $\recoursecontroltilde_{\openWclosed{\week}}$};
     \draw [ -latex,  dashed] (-5.25,1)-- (-5.25,0.15);
     \draw [ -latex,  dashed] (-5.25,1)-- (-4,0.15);
     \draw [ -latex,  dashed] (-5.25,1)-- (-2.25,0.15);
     \draw [ -latex,  dashed] (-5.25,1)-- (-1,0.15);
     \draw [ -latex,  dashed] (-5.25,1)-- (0.25,0.15);
     \draw [ -latex,  dashed] (-5.25,1)-- (2,0.15);
     \draw[dashed](-7,2.5) -- (-7,0);
     \draw[dashed](3.25,2.5) -- (3.25,1.4);
\end{tikzpicture}
}      \caption{Weekly decision-hazard-decision information structure.
    An arrow maps the available information towards the decision,
  so that, here, the information structure is partly anticipative as some arrows go from the right
  to the left.
  \label{Figure:WeeklyDHD}}
  \end{figure}

  Since the recourse decisions for the beginning of the week are made knowing the uncertainties for the whole week,
  they are anticipative.

For the weekly decision-hazard-decision structure, the information constraints are given by the $\sigma$-fields inclusions
  \begin{subequations}
  \begin{align}
    &\sigma\np{\va{\Controltilde}_{\closedWopen{\week}} } \subset
    \sigma\np{\va{\Uncertain}_{\winfhinf},
                \va{\Uncertain}_{\openWclosed{\binf{\week}}}, \dots,
                \va{\Uncertain}_{\openWclosed{\week\predecessor}} }
                \eqsepv
                \forall \week \in \WEEK
                \eqfinv
    \\
    &\sigma\np{\va{\RecourseControltilde}_{\openWclosed{\week}} } \subset
          \sigma\np{\va{\Uncertain}_{\winfhinf},
          \va{\Uncertain}_{\openWclosed{\binf{\week}}}, \dots,
          \va{\Uncertain}_{\openWclosed{\week\predecessor}}, \va{\Uncertain}_{\openWclosed{\week}} }
          \eqsepv
          \forall \week \in \WEEK          \eqfinp
  \end{align}
  \end{subequations}

  \subsection{Multistage stochastic optimization problem formulation}
\label{subsection:MultistageStochasticOptDHD}
  Now we classify the physical controls described in
 §\ref{subsection:variablesdefinition} into planning controls ${\controltilde}_{\closedWopen{\week}}$ and recourse controls
${\recoursecontroltilde}_{\openWclosed{\week}}$ as follows
\begin{subequations}
  \begin{align}
    \controltilde_{\closedWopen{\week}} &= {\techcommitmentSlow_{\closedWopen{\week}}}
                                          \eqfinv
  \\
  \recoursecontroltilde_{\openWclosed{\week}} &= \np{\techcommitmentFast_{\openWclosed{\week}},
                                        \tech_{\openWclosed{\week}},
                                      \pumping_{\openWclosed{\week}}, \turbining_{\openWclosed{\week}},
                                                \ens_{\openWclosed{\week}}}
                                                \eqfinp
  \end{align}
\end{subequations}

We model as planning decisions the on and off decisions for the ``slow'' thermal units, and as recourse decisions
all the remaining ones. We consider this classification to, in a way, model the rigidity of some thermal units,
 i.e. to model the fact that they cannot be switched on instantaneously.  Now the collection of on/off
 decisions for the thermal units is composed of $\techcommitment = \np{\techcommitmentSlow,\techcommitmentFast}$.

 In Sect.~\ref{section:AdequacyProblem}, we have presented the physical and economical model of the energy system.
 With the new notation in \S\ref{subsection:MultistageStochasticOptDHD}, we present in Table \ref{Table:DHDmathematicalnotation}
 the corresponding expressions and their compact mathematical versions.
  \begin{table}[!ht]
  \renewcommand{\arraystretch}{1.2}
  \centering
  \begin{tabular}{c|@{\hspace{0cm}}c@{\hspace{0.1cm}}|@{\hspace{0cm}}c@{\hspace{0.1cm}}|}
   \cline{2-3}
  &  \multicolumn{2}{c|}{Notation}\\
    \hline
    \multicolumn{1}{|c|}{              Function }                                                                                       &
                          Mathematical  &
                          Physical-Economical                                                                                                                     \\
                                                                                                                  \hline
                                                                                                                  \hline
  \multicolumn{1}{|c|}{
                        \begin{tabular}[c]{@{}c@{}}
                          \scriptsize        Storage \\
                          \scriptsize         dynamics
                    \end{tabular}
                            }                                         &
                      $\scriptstyle \widetilde{\Dynamics}_{{\week}}\bp{\stock_{\whinf}, \recoursecontroltilde_{\openWclosed{\week}}}$         &
                     $ \scriptstyle \dynamics_{\week}   \bp{ \stock_{\whinf}, \pumping_{\openWclosed{\week}} ,\turbining_{\openWclosed{\week}}}$
                                                                                                                             \\ \hline
  \multicolumn{1}{|c|}{
                        \begin{tabular}[c]{@{}c@{}}
                          \scriptsize         Thermal \\
                          \scriptsize            production
                    \end{tabular}
                            }                                         &
                    $\scriptstyle \widetilde{\Tech}_{\openWclosed{\week}}\bp{\controltilde_{\closedWopen{\week}}, \uncertain_{\openWclosed{\week}},\recoursecontroltilde_{\openWclosed{\week}}}$                       &
             $\scriptstyle \biggl\{\techproduction_{\wh\successor} \Bp{\bp{\techcommitmentSlow_{\wh}, \techcommitmentFast_{\wh\successor}}, \tech_{\wh\successor},\techuncertain_{\wh\successor}}\biggr\}_{\hour\in\HOUR} $
  \\ \hline
  \multicolumn{1}{|c|}{
                        \begin{tabular}[c]{@{}c@{}}
                  \scriptsize        Energy \\
                  \scriptsize       balance
                    \end{tabular}
                            }                                         &
                      $\scriptstyle \widetilde{\NodeBalance}_{\closedWopen{\week}}\bp{\controltilde_{\closedWopen{\week}}, \uncertain_{\openWclosed{\week}},\recoursecontroltilde_{\openWclosed{\week}}} $                       &
          \begin{tabular}[l]{@{\hspace{0cm}}p{6cm}@{\hspace{0cm}}}
              \raggedright $\scriptstyle  \biggl\{\nodeBalance \Bigl( \turbining_{\wh\successor},\pumping_{\wh\successor},  $
                  \\
              \centering   $\scriptstyle  \techproduction_{\wh\successor} \Bp{\bp{\techcommitmentSlow_{\wh}, \techcommitmentFast_{\wh\successor}}, \tech_{\wh\successor},\techuncertain_{\wh\successor}},$
                 \\
               \raggedleft $\scriptstyle   \uncertain^{\demand}_{\wh\successor},\ens_{\wh\successor}\Bigr)\biggr\}_{\hour\in\HOUR}$
          \end{tabular}
  \\ \hline
  \multicolumn{1}{|c|}{
                        \begin{tabular}[c]{@{}c@{}}
                          \scriptsize          Weekly \\
                          \scriptsize         cost
                    \end{tabular}
                            }
                   &
                  \hspace{0.01cm} $\scriptstyle \widetilde{\InstantaneousCost}_{\week}\bp{\controltilde_{\closedWopen{\week}},
                                                                                            \uncertain_{\openWclosed{\week}},
                                                                                            \recoursecontroltilde_{\openWclosed{\week}}}$
                  &
                  \begin{tabular}[l]{@{\hspace{0cm}}p{6cm}@{\hspace{0cm}}}
                    \raggedright  $ \scriptstyle \instantaneouscost_{\week} \Bp{
                           \bp{\techcommitmentSlow_{\closedWopen{\week}}, \techcommitmentFast_{\openWclosed{\week}} },
                       \tech_{\openWclosed{\week}}, \techuncertain_{\openWclosed{\week}}, \ens_{\openWclosed{\week}}}  = $
                 \\
                 \centering $
                 \scriptstyle\instantaneouscost \Bigl( \bp{\techcommitmentSlow_{\whinf},\techcommitmentFast_{\whinf\successor}},
                                                                                0,
                                                                                \techproduction_{\whinf\successor},
                                                                                \techuncertain_{\whinf\successor},
                                                                                 \ens_{\whinf\successor}
                \Bigr)+
           \displaystyle \sum_{\scriptscriptstyle \hour\in\HOUR \setminus\na{\binf{\hour}}}
           \scriptstyle\instantaneouscost \Bigl(\bp{\techcommitmentSlow_{\wh},\techcommitmentFast_{\wh\successor}},
           \bp{\techcommitmentSlow_{\wh\predecessor},\techcommitmentFast_{\wh}}     +
                                                                      $
             \\
              \raggedleft                      $
              \scriptstyle  \techproduction_{\wh\successor},\techuncertain_{\wh\successor},   \ens_{\wh\successor} \Bigr)$
                \end{tabular}
  \\ \hline
  \end{tabular}
  \caption{Correspondence between mathematical and physical and economical notations (Sect.~\ref{section:AdequacyProblem})
    in the weekly decision-hazard-decision framework 
  \label{Table:DHDmathematicalnotation}}
\end{table}
As already explained in~\S\ref{subsection:MultistageStochasticOptHD},
 we neglect the temporal coupling of thermal units from the last hour of one week to the first hour of the next week.
Considering the definitions in Table \ref{Table:DHDmathematicalnotation}, we formulate the physical
adequacy problem as a stochastic multistage optimization problem using the weekly decision-hazard-decision
information structure:
\begin{subequations}
  \label{eq:WeeklyPlanningWeeklyRecourseProblemFormulation}
  \begin{align}
    &
      \min_{ \va{\Stock}, \va{\Controltilde}, \va{\RecourseControltilde}}
        \EE  \biggl[
              \sum_{\week\in \WEEK}
              \widetilde{\InstantaneousCost}_{\week}\np{\va{\Stock}_{\whinf},
                                \va{\Controltilde}_{\closedWopen{\week}},
                                \va{\Uncertain}_{\openWclosed{\week}},
                                \va{\RecourseControltilde}_{\openWclosed{\week}}}
                  +
                  \FinalCost\np{\va{\Stock }_{\whlast}}			
                  \biggr]
                  \\
      & \suchthat \eqsepv
      \forall \week \in \WEEK
      \nonumber
      \\
      &
      \va{\Stock}_{\winfhinf} = \va{\Uncertain}_{\winfhinf} \eqfinv
      \\
      &
      \va{\Stock}_{\wsuccessorhinf} =\widetilde{\Dynamics}_{\week}\np{\va{\Stock}_{\whinf}, \va{\RecourseControltilde}_{\openWclosed{\week}}}
              \eqfinv
      \\
      &
      \widetilde{\NodeBalance}_{\closedWopen{\week}}\np{\va{\Controltilde}_{\closedWopen{\week}},
                                \va{\Uncertain}_{\openWclosed{\week}} ,  \va{\RecourseControltilde}_{\openWclosed{\week}}}
                                \geq
                                \va{0} \eqfinv \label{eq:balanceDHD}\\
      &\sigma\np{\va{\Controltilde}_{\closedWopen{\week}} } \subset
      \sigma\np{\va{\Uncertain}_{\winfhinf},
                  \va{\Uncertain}_{\openWclosed{\binf{\week}}}, \dots,
                  \va{\Uncertain}_{\openWclosed{\week\predecessor}} }
                  \eqfinv
      \\
      &\sigma\np{\va{\RecourseControltilde}_{\openWclosed{\week}} } \subset
            \sigma\np{\va{\Uncertain}_{\winfhinf},
            \va{\Uncertain}_{\openWclosed{\binf{\week}}}, \dots,
            \va{\Uncertain}_{\openWclosed{\week\predecessor}}, \va{\Uncertain}_{\openWclosed{\week}} }
            \eqfinp
   \end{align}
  \end{subequations}
  The final cost $\FinalCost\np{\va{\Stock }_{\whlast}}$ is used to give value to the energy in the storage at
  the end of the yearly period.

  \subsection{Bellman equations in decision-hazard-decision}
  \label{subsection:BellmanDHD}
  Defining the weekly state $\state_{\week} = \stock_{\whinf}$, that is, the storage level at the beginning of the week, we write the weekly Bellman equations~\eqref{eq:DHDweeklyBellman} using
 \cite[Proposition~13]{Carpentier-Chancelier-DeLara-Martin-Rigaut:2023}
\begin{subequations}
  \label{eq:DHDweeklyBellman}
  \begin{align}
    & \nBellman{\weeklast\successor}{\DHD}{\state_{\weeklast\successor}}=   \FinalCost\np{\state_{\weeklast\successor}}\eqfinv
    \\
    &\nBellman{\week}{\DHD}{\statew} =
    \min_{\controltilde_{\closedWopen{\week}}}
    {\EE} \bigg[
        {\min_{\recoursecontroltilde_{\openWclosed{\week}}} }
       \hspace*{0.1cm} \widetilde{\InstantaneousCost}_{\week}
        \np{\statew,
        \controltilde_{\closedWopen{\week}},
        \va{\Uncertain}_{\openWclosed{\week}},
        \recoursecontroltilde_{\openWclosed{\week}}}
        +
    \bBellman{\week\successor}{\DHD}{\widetilde{\Dynamics}_{\week}\np{
        \statew,\recoursecontroltilde_{\openWclosed{\week}}}}
                \bigg] \eqfinv \label{eq:DHDweeklyBellmanB}
\end{align}
\end{subequations}
 where the minimum inside the expectation term of Equation~\eqref{eq:DHDweeklyBellmanB}
 is computed subject to the constraint~\eqref{eq:balanceDHD}.

 If the sequence $\bp{\va{\Uncertain}_{\openWclosed{\binf{\week}}}, \dots,
\va{\Uncertain}_{\openWclosed{\week}} , \dots, \va{\Uncertain}_{\openWclosed{\bsup{\week}}}}$
of uncertainties is weekly independent,
the weekly Bellman equations provide an optimal solution for Problem~\eqref{eq:WeeklyPlanningWeeklyRecourseProblemFormulation}.
We highlight that, to get the optimal solution, the hourly uncertainties
  \(\va{\Uncertain}_{\openWclosed{\week}}=  \bp{\va{\Uncertain}_{\whinfsuccessor}, \dots, \va{\Uncertain}_{\wh\successor}, \dots, \va{\Uncertain}_{\wsuccessorhinf}} \)
within the week do not need to be assumed to be independent (from one hour to another).
Under this independence assumption,
it is well known that, for all~$\week\in \WEEK$, the function~$\nBellmanf{\week}{\DHD}$ satisfies
\begin{subequations}
  \begin{align}
  \nBellman{\week}{\DHD}{\statew} =
  &  \min_{ \substack{\va{\Controltilde}_{\closedWopen{\week}},\dots, \va{\Controltilde}_{\closedWopen{\bsup{\week}}}   \\
           \va{\RecourseControltilde}_{\openWclosed{\week}},\dots, \va{\RecourseControltilde}_{\openWclosed{\bsup{\week}}} }}
         \EE \biggl[
            \sum_{\week'=\week}^{\bsup{\week}}
            \InstantaneousCost_{\week'}\np{\va{\Stock}_{\wprimehinf},
            \va{\Controltilde}_{\closedWopen{\week'}},\va{\Uncertain}_{\openWclosed{\week'}},\va{\RecourseControltilde}_{\openWclosed{\week'}}
                        }
                +
                \FinalCost\np{\va{\Stock }_{\whlast}}			
                \biggr]
                \\
    & \suchthat \eqsepv
    \forall \week' \in \nc{\week,\dots,\bsup{\week}}
    \nonumber
    \\
    &
    \va{\Stock}_{\wprimehinf} =  \statew \eqfinv
    \\
    &
    \va{\Stock}_{\wprimesuccessorhinf} = \Dynamics_{\week'}\np{\va{\Stock}_{\wprimehinf},
                                                            \va{\RecourseControltilde}_{\openWclosed{\week'}}
                                                         }		
            \eqfinv
    \\
    &
    \NodeBalance_{\closedWopen{\week'}}\np{\va{\Controltilde}_{\closedWopen{\week'}},
    \va{\Uncertain}_{\openWclosed{\week'}}, \va{\RecourseControltilde}_{\openWclosed{\week'}}}
    \geq
         \va{0} \eqfinv \\
    &
    \sigma\np{ \va{\Controltilde}_{\closedWopen{\week'}}}
                        \subset \bsigmaf{
                                    \va{\Uncertain}_{\winfhinf},
                                    \va{\Uncertain}_{\openWclosed{\binf{\week}}},
                                    \dots,
                                    \va{\Uncertain}_{\openWclosed{\wprimepredecessor}}
                                    }
                  \eqfinv
                  \\
    &
    \sigma\np{ \va{\RecourseControltilde}_{\openWclosed{\week'}}}
                        \subset \bsigmaf{
                                    \va{\Uncertain}_{\winfhinf},
                                    \va{\Uncertain}_{\openWclosed{\binf{\week}}},
                                    \dots,
                                    \va{\Uncertain}_{\openWclosed{\wprimepredecessor}},
                                    \va{\Uncertain}_{\openWclosed{\week'}}
                                    }
                  \eqfinv
 \end{align}
\end{subequations}
so that, the  value $ \nBellman{\week}{\DHD}{\statew}$ of Bellman function is interpreted as the future optimal cost
when at week $\week$ the storage level is $\statew$ and the weekly decision-hazard-decision
information structure is considered.


\subsection{Theoretical comparison between Bellman functions in $\HD$ and $\DHD$}
\label{subsection:TheoreticalBellmanComparison}
 When we compare theoretically the Bellman functions $\nBellmanf{\week}{\HD}$ and $\nBellmanf{\week}{\DHD}$, given
 by the Bellman equations~\eqref{eq:HDweeklyBellman} and~\eqref{eq:DHDweeklyBellman},
 we observe that the weekly hazard-decision  approach is a relaxation of the weekly decision-hazard-decision  approach
  with respect to the information constraint.
In other words, at each stage of the stochastic multistage optimization problem, the decision maker has more
information when making the decision in the~$\HD$ case, than in the~$\DHD$ case.
Therefore, we have the following inequalities for all~$\week\in\WEEK$: 
\begin{align}
  \label{eq:BellmanRelation}
  \nBellmanf{\week}{\HD} \leq \nBellmanf{\week}{\DHD}
  \eqfinp
\end{align}

The usage values --- or prices~$\usagevalue$ --- are defined as the opposite of the derivative of the Bellman functions with
respect to the storage level \cite{Steeger2014}:
\begin{subequations}
  \label{eq:usagevaluesdefinition}
\begin{align}
  \usagevalue^{\HD}_{\week} &= - \frac{\mathrm{d}}{\mathrm{d}\state} \nBellman{\week}{\HD}{\state} \eqfinv
\\
  \usagevalue^{\DHD}_{\week} &= - \frac{\mathrm{d}}{\mathrm{d}\state} \nBellman{\week}{\DHD}{\state} \eqfinp
\end{align}
\end{subequations}
Whereas we have established an inequality between Bellman functions in Equation~\eqref{eq:BellmanRelation},
it is impossible to do so for the usage values (as they are derivatives).

In Sect.~\ref{section:numericalstudy}, we compare numerically the effect of
both information structures when computing usage values for a case study of small size.

\section{Numerical study}
\label{section:numericalstudy}

In this Sect.~\ref{section:numericalstudy}, we present a numerical study. In \S\ref{subsection:StudyCaseDescription}, we describe
the energy system of small size considered. In \S\ref{subsection:BellmanComputation},
we present the method to compute the Bellman functions and usage values, either for the hazard-decision information structure in
Sect.~\ref{section:AdequacyProblemHD} or for the decision-hazard-decision information structure in Sect.~\ref{section:AdequacyProblemDHD}.
In \S\ref{subsection:NumericalResultsUsageValues}, we compare the numerical results for the usage values.
Finally, in \S\ref{subsection:Simulation}, we show the effect of the information modelling choice when designing policies
to carry on a simulation of the energy system dispatch.

\subsection{Case study description}
\label{subsection:StudyCaseDescription}

We consider a small energy system to conduct the numerical study with the following components:
3~thermal units (base unit, semi-base unit and peak unit);
1~residual demand;
1~storage unit.
The thermal units differ in their prices and power capacities. Whereas the base and semi-base units are the cheapest
and the less flexible ones, the peak unit is the most expensive but the most flexible one.

We model the components mentioned above as in
Sect.~\ref{section:AdequacyProblem}, and we consider  the
\emph{weekly hazard-decision} (Sect.~\ref{section:AdequacyProblemHD}) and the
\emph{weekly decision-hazard-decision}  (Sect.~\ref{section:AdequacyProblemDHD}) mathematical modelling of the problem.
In the weekly hazard-decision case, all the decisions for the week are anticipative but, in the
weekly decision-hazard-decision case, we consider that
the on and off decisions for the base and semi-base units are nonanticipative (planned in advance).
This classification is made taking into account the rigidity and flexibility of the units.

To model the uncertainties, we have taken a finite number of scenarios
 provided by RTE, with uniform probability.
We divide the scenarios into
two sets:
  \begin{itemize}
    \item  $\bseqp{\uncertain_{\openWclosed{\week}}^n}{\week\in\WEEK}$ for $n\in\ic{1,N}=\na{1,\ldots,N}$ used to compute Bellman functions and their associated
    nonanticipative policies,
    \item  $\bseqp{\uncertain_{\openWclosed{\week}}^\chronicle}{\week\in\WEEK}$ for $\chronicle\in\ic{1,C}=\na{1,\ldots,C}$ used in simulation.
  \end{itemize}

For the computations, we use 
  the programming language Julia \cite{bezanson2017julia}, the  JuMP package~\cite{Lubin2023} and the Xpress solver.

\subsection{Bellman functions computation}
\label{subsection:BellmanComputation}

We compute the Bellman functions using the classical stochastic dynamic programming
(SDP) algorithm \cite{Bertsekas:2000}
for the hazard-decision framework, and an adapted version for
the decision-hazard-decision framework.
In both cases, the problem for the entire year is decomposed
in weekly problems that we solve using the backward recursions~\eqref{eq:HDweeklyBellman} and
\eqref{eq:DHDweeklyBellman} (respectively) for a regular discretization grid~$\GRID_\STATE$
of the state space~$\STATE=\nc{\binf{\stock}, \bsup{\stock}}$ (see Table~\ref{Table:storagevariables}).
This discretization grid has been refined to achieve a satisfactory piecewise approximation of the Bellman functions.

On the one hand, we do it to be able to handle
binary or discrete variables within the planning decisions, in which case SDDP is not suitable
(but variants like SDDiP are).
On the other hand, we keep the original state without augmenting it as in ASDH-SDDP in~\cite{Street-Valladao-Lawson-Velloso:2020}.
The expectations in the Bellman equations  are computed as finite sums.

The objective is to obtain  two sequences \(\bseqp{\nBellmanf{\week}{\HD}}{\week\in\WEEK}\) and
\(\bseqp{\nBellmanf{\week}{\DHD}}{\week\in\WEEK}\) of Bellman functions corresponding to the
weekly hazard-decision and weekly decision-hazard-decision frameworks. They are represented by real numbers
\(\bseqp{\nBellman{\week}{\HD}{\state}}{\week\in\WEEK}\) and
\(\bseqp{\nBellman{\week}{\DHD}{\state}}{\week\in\WEEK}\) with $\state$ varying in a grid~$\GRID_\STATE$.
For both information structures, solving the Bellman equations leads to MILP problems for each week.
Indeed, we approximate the cost-to-go functions --- $\nBellmanf{\week\successor}{\HD}$ and
  $\nBellmanf{\week\successor}{\DHD}$
  in~\eqref{eq:HDweeklyBellman}  and~\eqref{eq:DHDweeklyBellman} respectively ---
  as continuous piecewise linear functions,
  but not necessarily convex (due to the presence of binary decisions).
  Taking into account piecewise linear cost-to-go functions in MILPs can be handled
  through the so-called lambda method \cite{LEE2001269}.
The algorithms are sketched in Alg.~\ref{alg:HD} and  Alg.~\ref{alg:DHD}.
\begin{algorithm}[ht]
\caption{Hazard-decision stochastic dynamic programming \label{alg:HD}}
\KwData{uncertainties scenarios $\uncertain_{\openWclosed{\week}}^n$ with $n \in  \ic{1, N}$ and probabilities $\frac{1}{N} $, \\
       \hspace{0.94cm}space state discretization $\state \in \GRID_\STATE$,  \\
        \hspace{0.94cm}final cost $\FinalCost\bp{\state_{\weeklastsuccessor}}$ }
 \KwResult{\(\bseqa{\nBellman{\week}{\HD}{\state}}{\week\in\WEEK}\)}
 \For{$\week = \bsup{\week},  \dots, \binf{\week}$}{
    \For{$\state \in \GRID_\STATE$}{
      $\nBellman{\week}{\HD,0}{\state} = 0$; \\
        \For{$n \in \ic{1, N}$}{
          $\nBellman{\week}{\HD,n}{\state} = \nBellman{\week}{\HD,n-1 }{\state}+
          \frac{1}{N}
           \displaystyle  \min_{\recoursecontrol_{\openWclosed{\week}}^n}
          \Bigl\{\InstantaneousCost_{\week}\bp{
                                              \state,
                                            \uncertain_{\openWclosed{\week}}^n,
                                            \recoursecontrol_{\openWclosed{\week}}^n}
            +
              \bBellman{\week\successor}{\HD}{\Dynamics_{\week}\np{\state,\recoursecontrol_{\openWclosed{\week}}^n}}\eqsepv $ \\
              \hspace{7cm} $\text{such that ~\eqref{eq:balanceHD}}
          \Bigr\}$
     }
  } }
\end{algorithm}

  \begin{algorithm}[ht]
    \caption{Decision-hazard-decision stochastic dynamic programming}\label{alg:DHD}
    \KwData{uncertainties scenarios $\uncertain_{\openWclosed{\week}}^n$ with $n \in \ic{1, N}$ and probabilities $\frac{1}{N} $, \\
         \hspace{0.94cm}space state discretization $\state \in \GRID_\STATE$,  \\
          \hspace{0.94cm}final cost $\FinalCost\bp{\state_{\weeklastsuccessor}}$ }
   \KwResult{\(\bseqa{\nBellman{\week}{\DHD}{\state}}{\week\in\WEEK}\)}
   \For{$\week = \bsup{\week},  \dots, \binf{\week}$}{
      \For{$\state \in \GRID_\STATE$}{
            $\nBellman{\week}{\DHD}{\state} = \displaystyle  \min_{\controltilde_{\closedWopen{\week}}}  \sum_{n=1}^N \frac{1}{N}
           \biggl(  \min_{\recoursecontroltilde_{\openWclosed{\week}}^n}
            \Bigl\{\InstantaneousCost_{\week}\bp{
                \state,
                \controltilde_{\closedWopen{\week}},
                \uncertain_{\openWclosed{\week}}^n,
               \recoursecontroltilde_{\openWclosed{\week}}^n} + \bBellman{\week\successor}{\DHD}{\Dynamics_{\week}\np{\state,\recoursecontroltilde_{\openWclosed{\week}}^n} }\eqsepv $\\
               \hspace{8.1cm} $\text{such that~\eqref{eq:balanceDHD}}
                \Bigr\}
                  \biggr)$
       }
    }
    \end{algorithm}

    In the Alg.~\ref{alg:HD} ($\HD$), as many deterministic MILPs  are solved
    as there are scenarios, at each week and each point of the grid~$\GRID_\STATE$.
In the Alg.~\ref{alg:DHD} ($\DHD$), a stochastic two-stage MILP is solved
using its extensive formulation.

When comparing the Bellman functions for all the weeks in the year, we obtain that
$\nBellman{\week}{\HD}{\state} < \nBellman{\week}{\DHD}{\state}$ for all $\state$ in the grid~$\GRID_\STATE$, which
  numerically confirms the theoretical result in ~\eqref{eq:BellmanRelation}.

The usage values or prices, $\usagevalue$, are calculated as the opposite of the derivative of the Bellman functions with
respect to the storage level as in Equation~\eqref{eq:usagevaluesdefinition}.
As the Bellman functions are computed on the discrete grid~$\GRID_\STATE$,
and as usage values in~\eqref{eq:usagevaluesdefinition} have been defined as
derivatives,
usage values are approximated by increments between \emph{middle points~$\state$} of
 the state space grid~$\GRID_\STATE$ by the formula
\begin{equation}
  \frac{\mathrm{d}}{\mathrm{d}\state} \nBellman{\week}{}{\state}  \approx \frac{\nBellman{\week}{}{\state + \frac{\Delta \state}{2}}  -
                                                     \nBellman{\week}{}{\state-\frac{\Delta \state}{2}}}
                                                     {\Delta\state}
                                                     \eqfinv
\end{equation}
with $\Delta \state$ the discretization step of the grid~$\GRID_\STATE$.

\subsection{Numerical results for usage values}
\label{subsection:NumericalResultsUsageValues}

\begin{figure}[htbp]
\includegraphics[width=0.9\columnwidth]{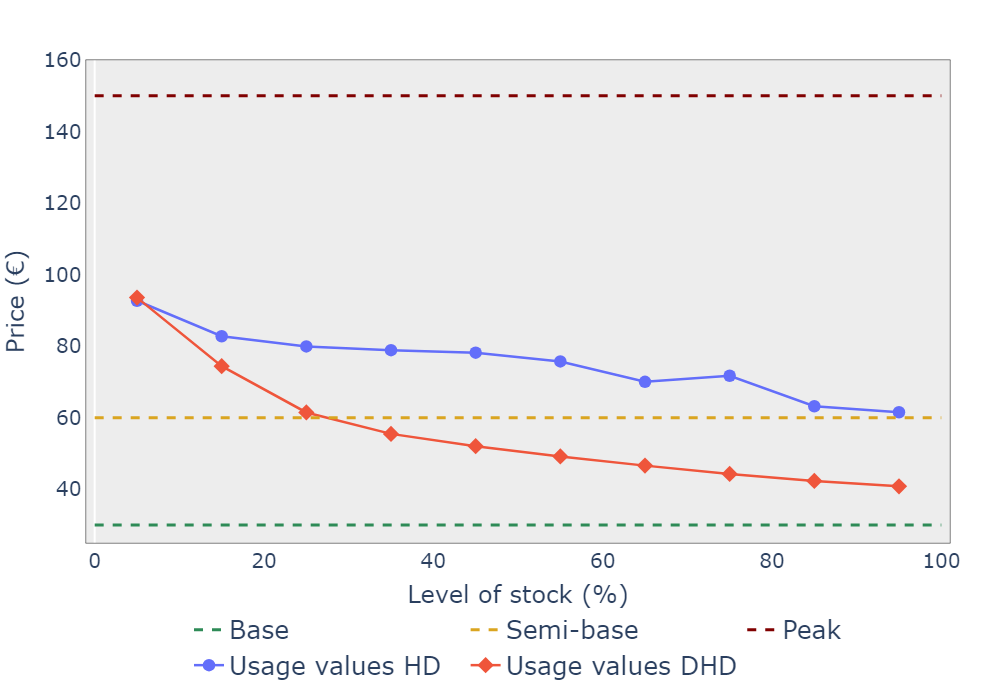}
\caption{Prices and usage values comparison for week~20,
  depending on the information structure ($\DHD$ or $\HD$) considered in the Bellman equations
\label{Figure:PricesWeek2}}
\end{figure}

In Fig.~\ref{Figure:PricesWeek2}, we compare the prices of the
thermal units with the usage values obtained with the Bellman functions computed using the
$\HD$ and $\DHD$ structures.

  We observe that, for storage levels higher than 20\%, the $\HD$ usage value in blue is above
    the semi-base unit price (dashed yellow), whereas the $\DHD$ usage value in red is below.
  As a consequence, we expect different dispatches in simulation when
  using these usage values to design the storage policy, since the merit order of the production means changes.
  This change in the dispatch could lead to different conclusions when carrying out prospective studies.
  In \S\ref{subsection:Simulation} we compare the results in simulation
  induced by both modelling options to compute Bellman functions.

  \subsection{Numerical comparison between policies in $\HD$ and $\DHD$}
  \label{subsection:Simulation}
  The goal of this part of the numerical study is to compare
  the dispatches obtained when using the policies induced by both Bellman functions
  \(\bseqp{\nBellmanf{\week}{\HD} }{\week\in\WEEK}\) and
  \(\bseqp{\nBellmanf{\week}{\DHD} }{\week\in\WEEK}\)
  computed in
  §\ref{subsection:BellmanComputation} .

  We recall that an uncertainty scenario for simulation (also called \emph{chronicle})  is
  denoted by \(\bseqp{\uncertain_{\openWclosed{\week}}^{\chronicle}}{\week\in\WEEK} \).
  We simulate the operation of the system
  for different chronicles $\chronicle \in\ic{1,C}$ in which the availability and demand change.

  For the sake of simplicity and to be consistent with the weekly decision-hazard-decision
  information structure, we choose to design an ``anticipative simulator'': uncertainties are disclosed
  at the end of the week and recourse controls are computed at the end of the week but applied within the week.
  The simulation algorithm is the same regardless the sequence of Bellman functions ($\HD$ or $\DHD$) chosen to design the policies.
  We illustrate how the
  simulation is done for one uncertain chronicle in  Algorithm~\ref{alg:Simulation}.
  We highlight that this simulation algorithm can be used with any Bellman functions.
  \begin{algorithm}[ht!]
    \caption{Decision-hazard-decision policy and simulation}\label{alg:Simulation}
  \KwData{uncertainties scenarios $\uncertain_{\openWclosed{\week}}^n$ with $n \in \ic{1, N}$ and probabilities $\frac{1}{N} $, \\
          \hspace{0.94cm}Bellman functions $\bseqa{\nBellmanf{\week}{}}{\week\in\WEEK}$,\\
          \hspace{0.94cm}simulation chronicle $\bseqa{\uncertain_{\openWclosed{\week}}^{\chronicle}}{\week\in\WEEK}$,\\
          \hspace{0.94cm}initial condition for the state $\state_{0}$}
   \KwResult{\(\Bseqa{\statew^{\chronicle}, \control_{\closedWopen{\week}}^{\chronicle},  \recoursecontrol_{\openWclosed{\week}}^{\chronicle} }{\week\in\WEEK}\)}
   $\state_{\binf{\week}}^{\chronicle} = \state_{0}$\\
   \For{$\week = \binf{\week},  \dots, \bsup {\week}$}{
    \texttt{compute nonanticipative controls:}\\
      \(\displaystyle \control_{\closedWopen{\week}}^{\chronicle} =  \argmin_{\controltilde_{\closedWopen{\week}}}  \sum_{n=1}^N \frac{1}{N}
      \bgp{  \min_{\recoursecontroltilde_{\openWclosed{\week}}^n}
       \Ba{\InstantaneousCost_{\week}\bp{
           \statew^{\chronicle},
           \controltilde_{\closedWopen{\week}},
           \uncertain_{\openWclosed{\week}}^n,
          \recoursecontroltilde_{\openWclosed{\week}}^n} + \bBellman{\week\successor}{}{\Dynamics_{\week}\np{\statew^{\chronicle},\recoursecontroltilde_{\openWclosed{\week}}^n} }
    }
    }\)
       \\
       \texttt{compute recourse controls:}\\
   \(\displaystyle \recoursecontrol_{\openWclosed{\week}}^{\chronicle} =  \argmin_{\recoursecontroltilde_{\openWclosed{\week}}}
        \Ba{\InstantaneousCost_{\week}\bp{
            \statew^{\chronicle},
            \control_{\closedWopen{\week}}^{\chronicle},
            \uncertain_{\openWclosed{\week}}^{\chronicle},
            \recoursecontroltilde_{\openWclosed{\week}}} + \bBellman{\week\successor}{}{\Dynamics_{\week}\np{\statew^{\chronicle}
                                                                                                      ,\recoursecontroltilde_{\openWclosed{\week}}} }
    }\)
    \\
    \texttt{update state:}\\
    \(\state_{\week\successor}^{\chronicle} = \Dynamics_{\week}\np{\statew^{\chronicle},\recoursecontroltilde_{\openWclosed{\week}}}\)
       }
    \end{algorithm}


  \begin{figure}[htbp]
            \begin{subfigure}[b]{\columnwidth}
                  \includegraphics[width=0.9\columnwidth]{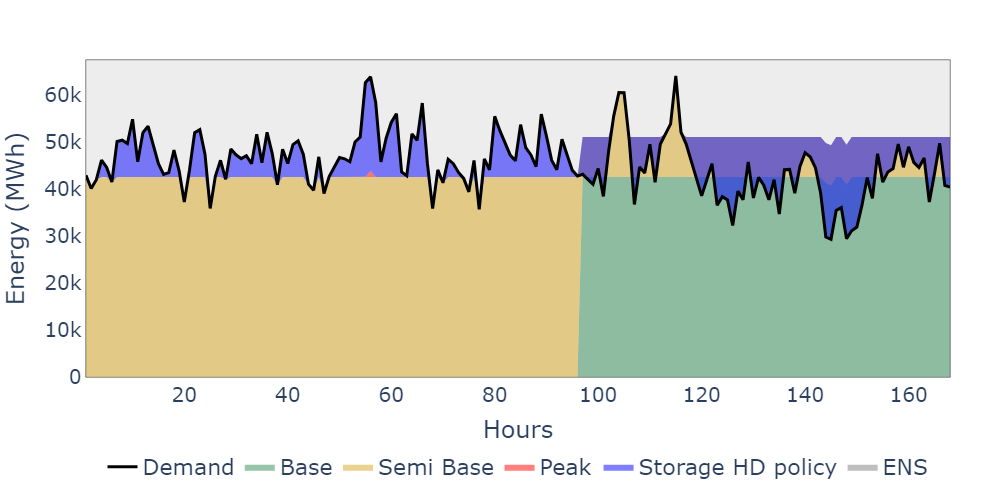}
                  \caption{Weekly hazard-decision policy dispatch
                  \label{fig:subfigHD}}
            \end{subfigure}
            \begin{subfigure}[b]{\columnwidth}
                \includegraphics[width=0.9\columnwidth]{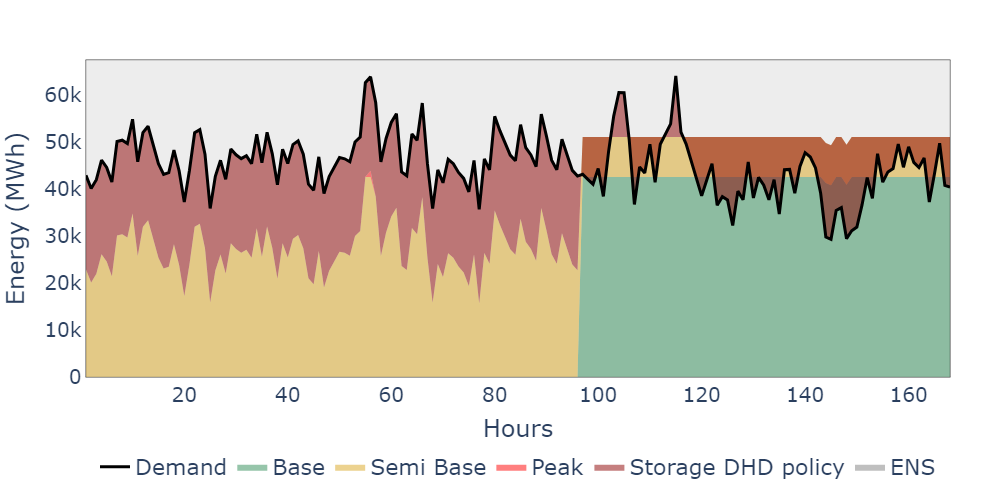}
                \caption{Weekly decision-hazard-decision policy dispatch
                \label{fig:subfigDHD}}
            \end{subfigure}
            \caption{Dispatch comparison using the policies induced by the Bellman functions computed
              either with the weekly hazard-decision (Fig.~\ref{fig:subfigHD})
or with the weekly decision-hazard-decision information structure (Fig.~\ref{fig:subfigDHD})
     \label{Figure:DispatchWeek2}}
    \end{figure}

    In  Fig.~\ref{Figure:DispatchWeek2}, we compare the simulations of the dispatches
    for one chronicle of uncertainties (residual demand and availability fixed) in the week 20 of the year when
     using the $\HD$ and $\DHD$ Bellman functions to compute the policies. In the selected chronicle, the base unit
       is not available until the hour 96 of the week. The semi-base and peak units are available during the entire week.

    We compare the dispatch obtained with the weekly hazard-decision policy in Fig.~\ref{fig:subfigHD} with the one
    obtained with the weekly decision-hazard-decision policy in Fig.~\ref{fig:subfigDHD}.
    To satisfy the same demand, with the same production means available, we observe that the dispatch varies depending
     on the policy considered.     In the case of the $\HD$ policy, the semi-base unit works at its maximum power and
     the marginal production is made by the storage.
     On the contrary, in the $\DHD$ policy dispatch, the storage operates at its maximum power and the marginal
     production is made by the semi-base unit.

    This difference comes from the fact that, as we remarked in Fig.~\ref{Figure:PricesWeek2},
    the $\HD$ usage values are above the semi-base unit price,
     whereas the {$\DHD$} usage values are below the semi-base unit price. This changes the dispatch order.

     The merit order difference in the two policies studied is observed for several  weeks in the year, which leads,
     for all chronicles, to a higher use of the storage (more turbining and more pumping) when using the {$\DHD$} policy.

  \begin{figure}[htbp]
      \includegraphics[width=0.9\columnwidth]{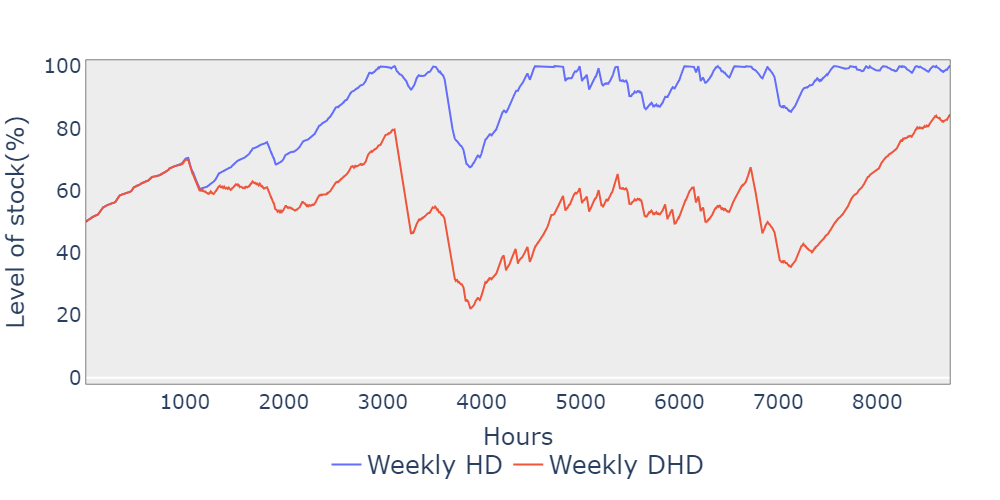}
      \caption{Comparison of the storage level obtained from the weekly hazard-decision and the
      weekly decision-hazard-decision policies for one given chronicle (of uncertainties)
      \label{Figure:LevelStockYear}}
    \end{figure}

     However, we observe that the $\HD$ policy tends to store more energy than the $\DHD$ policy as
     we can see in Fig.~\ref{Figure:LevelStockYear} for the same simulated yearly chronicle. This result is observed in all
    simulated chronicles and can be explained with the usage value   difference: in the case of the $\DHD$~policy, there
     is less  interest in storing energy since it has lower usage value.


\section{Conclusion}
\label{Section:Conclusion}

In this paper, we have first formalized the weekly hazard-decision
information structure in a two-timescale setting.
The weekly hazard-decision information structure is the current reference framework when modelling
the resource adequacy problem as a multistage stochastic optimization problem to compute usage values for prospective
studies  under uncertainty. We have written the correponding weekly Bellman equations, that make it possible to compute usage
 values for the energy in the storage. These Bellman equations respect the hourly physical constraints.

Then, we have highlighted the need to improve this structure to account for temporal rigidities in thermal operation, that is,
  the fact that on/off decisions cannot be modelled as fast or last minute decisions.
  Therefore, we have introduced the weekly decision-hazard-decision information structure
  in a two-timescale setting.
In this structure, the decisions for each stage in the decision-making process are separated
 into planning and recourse decisions, depending on the physical modelling aspects of the decisions.
We have presented the mathematical formulation of the problem considering the weekly decision-hazard-decision information
structure  and its associated Bellman equation in the weekly timescale,
 that still respect the hourly physical constraints.

 Afterwards, we have have carried out a numerical analysis to quantify, in a case study, the consequences of the
 information structure modelling choice when computing Bellman functions and their associated usage values.
 From this study, we have inferred that the policy induced by the weekly hazard-decision information modelling leads to an
 overestimation  of the thermal flexibility and as a result, a lower use of the storage.
 Indeed, when using the policy induced by the weekly decision-hazard-decision Bellman functions,
 we take into account some rigidity in the thermal units (specially  on/off decisions) and the storage is more used.
In consequence, we concur with~\cite{Street-Valladao-Lawson-Velloso:2020} that decision-hazard-decision information structures are of great interest to calculate
usage values, taking into account that some types of thermal units are less flexible than others.

This study has been carried using a single storage, so that the well-known curse of dimensionality in SDP is not binding.
In future work, we will turn to spatial decomposition techniques \cite{Pacaud-DeLara-Chancelier-Carpentier:2022}
to extend the study to multiple storage facilities.

\newcommand{\noopsort}[1]{} \ifx\undefined\allcaps\def\allcaps#1{#1}\fi


\appendix  
\section{Hourly composition of the storage dynamics}
\label{Appendix:Hourlycomposition}
We consider the hourly dynamics function given by
 \begin{subequations}
  \begin{align}
    \label{eq:balanceExp}
\dynamics\np{ \stock, \pumping ,\turbining
}
= \stock + \eta \pumping - \turbining
\eqfinp
\end{align}

This equation represents the hourly evolution of the level in a storage.
 To be able to formulate the problem in a weekly framework, we define the weekly dynamics
 $\dynamics_{\week}$.
 The weekly dynamics  $\dynamics_{\week}\np{ \stock_{\whinf}, \pumping_{\openWclosed{\week}},\turbining_{\openWclosed{\week}}}$
 is obtained by composition of the hourly  dynamics~$\dynamics\np{ \stock_{\wh}, \pumping_{\wh\successor} ,\turbining_{\wh\successor}} $
 for $\hour \in \HOUR$, but the composition being done only on storage level variable.

 To describe the dynamics composition we introduce the following extra notation
\begin{align}
    \dynamics^{\pumping, \turbining}_{\wh}
                            \bp{\stock_{\wh}}
        =
            \dynamics
                            \bp{\stock_{\wh}, \pumping_{\wh\successor}, \turbining_{\wh\successor}}
    \eqfinv
\end{align}
so that we get
\begin{align}
    \dynamics_{\week}\bp{\stock_{\whinf}, \pumping_{\openWclosed{\week}},\turbining_{\openWclosed{\week}}
                            }
        =\Bp{
            \dynamics^{\pumping, \turbining}_{\couple{\week}{\bsup{\hour}}}
                \circ
                    \dynamics^{\pumping, \turbining}_{\couple{\week}{\bsup{\hour}\predecessor}}
                \circ
                    \cdots
                \circ
                    \dynamics^{\pumping, \turbining}_{\whinf}
                    }
                            \np{\stock_{\whinf}
                                }
    \eqfinp
\end{align}
From the hourly dynamics~$\dynamics$ expression~\eqref{eq:balanceExp}, the weekly dynamics $\dynamics_{\week}$ is given by:
\begin{align}
  \dynamics_{\week}\bp{\stock_{\whinf}, \pumping_{\openWclosed{\week}},\turbining_{\openWclosed{\week}}
  } =
  \stock_{\whinf} + \sum_{\hour\in\HOUR} \bp{\eta \pumping_{\wh\successor} - \turbining_{\wh\successor}} \eqfinp
\end{align}

 \end{subequations}

\section{Ideal information structure in decision-hazard-decision}
\label{Appendix:HourlyRecourseBellman}

As discussed in~\S\ref{Our_contribution}, we should ideally
consider an information structure in which, at the beginning of the week~$\week$,
the planning for the
nonanticipative decisions~$\controltilde_{\closedWopen{\week}}$ is made knowing
only the past uncertainties, that is, the uncertainties up to
$\whinf$. Then, when the uncertainties
$\uncertain_{\wh}$ begin to be disclosed hour by hour, there are so-called
hourly recourse decisions~$\recoursecontroltilde_{\wh}$
that are also made hour by hour knowing the uncertainties up to $\wh$. Such information structure, illustrated in Fig.~\ref{Figure:WeeklyDHourlyHD}, is
called \emph{decision-hazard-decision information structure with hourly recourse}.

\begin{figure}[htbp]
  \centering
  \resizebox{0.8\columnwidth}{!}{
    \begin{tikzpicture}
     \filldraw[black] (-7.,0) circle (2pt) node[anchor=north] {$\winfhinf$};
     \draw[dashed, thick] (-7.,0) -- (-5.25,0);
     \node (us) at (-5.25,1.2) {\footnotesize {${\controltilde_{\closedWopen{\week}}}$}};
     \draw [decorate,
           decoration = {brace, raise=5pt,
               amplitude=8pt}] (-7.,1.25) --  (-5.25,1.25);
           \node (W-DH) [align=left] at (-5.5,2.2) { \scriptsize Information before  \\ \scriptsize making \scriptsize decision $\controltilde_{\closedWopen{\week}}$};
     \filldraw[black] (-5.25,0) circle (2pt) node[anchor=north] {$\whinf$};
     \draw[black, thick] (-5.25,0) -- (-4,0);
     \filldraw[black] (-4,0) circle (2pt) node[anchor=north] {$\whinfsuccessor$};
     \draw[dashed, thick] (-4,0) -- (-2.25,0);
     \filldraw[black] (-2.25,0) circle (2pt) node[anchor=north] {$\wh\predecessor$};
     \draw[black, thick] (-2.25,0) -- (-1,0);
     \filldraw[black] (-1,0) circle (2pt) node[anchor=north] {$\wh$};
     \draw[black, thick] (-1,0) -- (0.25,0);
     \filldraw[black] (0.25,0) circle (2pt) node[anchor=north] {$\wh\successor$};
     \draw[dashed, thick] (0.25,0) -- (2,0);
     \filldraw[black] (2,0) circle (2pt) node[anchor=north] {$\whsup$};
     \draw[dashed, thick] (0.25,0) -- (2,0);
     \draw[black, thick] (02 ,0) -- (3.25,0);
     \filldraw[black] (3.25,0) circle (2pt) node[anchor=north] {$\wsuccessorhinf$};
     \draw[dashed, thick] (3.25,0) -- (4.5,0);
 \draw [decorate,
 decoration = {brace, raise=5pt, mirror,
     amplitude=8pt}] (-4,-0.5) --  (0.25,-0.5);
 \node (uncertain)[align=center] at (-1.875, -1.2) {\footnotesize	 Uncertainties};
 \node () at (-1.875, -1.5){\footnotesize	 $\bp{\uncertain_{\whinfsuccessor}, \dots, \uncertain_{\wh\predecessor}, \uncertain_{\wh}, \uncertain_{\wh\successor}}$};
          \draw[ -latex,  dashed] (-4 ,1)-- (-4, 0.15);
          \node (1) at (-4,1.2) {\footnotesize {${\recoursecontroltilde_{\whinfsuccessor}}$}};
          \draw[ -latex,  dashed] (-2.25 ,1)-- (-2.25,0.15);
          \node (1) at (-2.25,1.2) {\footnotesize {${\recoursecontroltilde_{\wh\predecessor}}$}};
          \draw[ -latex,  dashed] (-1 ,1)-- (-1,0.15);
          \node (1) at (-1,1.2) {\footnotesize {${\recoursecontroltilde_{\wh}}$}};
          \draw[ -latex,  dashed] (0.25 ,1)-- (0.25,0.15);
          \node (1) at (0.25,1.2) {\footnotesize {${\recoursecontroltilde_{\wh\successor}}$}};
 \draw [decorate,
        decoration = {brace, raise=5pt,
         amplitude=8pt}] (-7,2.5) --  (0.25,2.5);
     \node (W-HD)[align=left] at (-3.4,3.2) { \scriptsize Information before making recourse decision $\recoursecontroltilde_{\wh\successor}$};
     \draw [ -latex,  dashed] (-5.25,1)-- (-5.25,0.15);
     \draw [ -latex,  dashed] (-5.25,1)-- (-4,0.15);
     \draw [ -latex,  dashed] (-5.25,1)-- (-2.25,0.15);
     \draw [ -latex,  dashed] (-5.25,1)-- (-1,0.15);
     \draw [ -latex,  dashed] (-5.25,1)-- (0.25,0.15);
     \draw [ -latex,  dashed] (-5.25,1)-- (2,0.15);
     \draw[dashed](-7,2.5) -- (-7,0);
     \draw[dashed](0.25,2.5) -- (0.25,1.4);
\end{tikzpicture}
}
\caption{Decision-hazard-decision information structure with hourly recourse.
  An arrow maps the available information towards the decision,
  so that, here, the information structure is nonanticipative as all arrows go either down or from the left to the right.
\label{Figure:WeeklyDHourlyHD}}
\end{figure}

This structure leads to Bellman equations (see
\S\ref{Bellman_equations_in_decision-hazard-decision_information_structure_with_hourly_recourse})
with nested mathematical expectations and minimizations, hence out of reach
to compute. This is why we opted in Sect.~\ref{section:AdequacyProblemDHD}
for a compromise solution, namely, the weekly decision-hazard information structure.
\begin{remark}
  In real operation it is possible to adjust the production ``in real time" to meet demand, that is,
  there is some anticipativity at the fast timescale for the hourly recourse decision~$\recoursecontroltilde_{\wh}$. Therefore,
it is not unrealistic to consider the ideal model described above.
Note that, it is not possible to consider a fully nonanticipative (or decision-hazard)  structure,
in which  the hourly recourse decisions  $\recoursecontroltilde_{\wh}$ are made knowing the uncertainties up to
$\wh\predecessor$ (instead of $\wh$),
 because it would lead to an excessive usage of the slack variable $\ens$ in the energy balance equation~\eqref{eq:BalanceEquations}.
\end{remark}

\subsection{Weekly decision-hazard-decision information structure with hourly recourse}

We consider an information structure in which  the nonanticipative control~$\controltilde_{\closedWopen{\week}}$
is made knowing the uncertainties up to the beginning of week~$\week$,
and the hourly recourse control~$\recoursecontroltilde_{\wh\successor}$ is made knowing  the uncertainties
up to the hour~$\wh\successor$.
As illustrated in Fig.~\ref{Figure:WeeklyDHourlyHD},
at the beginning of the week~$\week$ the nonanticipative decision
$\controltilde_{\closedWopen{\week}}$ is made (for all $\hour\in\HOUR$) and then,
once the uncertainties are disclosed \emph{hour by hour}, the recourse decisions
$\recoursecontrol_{\wh\successor}$ are sequentially made.

For the weekly decision-hazard-decision structure  with hourly recourse the information constraints
are given by the $\sigma$-fields inclusions
\begin{subequations}
\begin{align}
  &\forall  \week\in\WEEK \eqfinv  \nonumber \\
  &&    \bsigmaf{\va{\Controltilde}_{\closedWopen{\week}}}
          &\subset
             \bsigmaf{
                \va{\Uncertain}_{\winfhinf},
                \va{\Uncertain}_{\openWclosed{\binf{\week}}},
                \dots,
                \va{\Uncertain}_{\openWclosed{\week\predecessor}}
             }
             \eqsepv
\\
&\forall \wh\in\timeline  \eqfinv \nonumber \\
 &&  \bsigmaf{\va{\RecourseControltilde}_{\wh\successor}}
        &\subset
            \bsigmaf{
                \va{\Uncertain}_{\winfhinf},
                \va{\Uncertain}_{\openWclosed{\binf{\week}}},
                \dots,
                \va{\Uncertain}_{\openWclosed{\week\predecessor}},
                \va{\Uncertain}_{\whinfsuccessor},
                \dots,
                \va{\Uncertain}_{\wh\successor}
                }
                \eqfinp
    \label{eq:Info-WPHR-V}
\end{align}
\end{subequations}

\subsection{Multistage stochastic optimization problem formulation}

Considering the definitions in Table~\ref{Table:DHDmathematicalnotation}, we formulate the physical-economical
adequacy problem as a stochastic multistage optimization problem using the weekly decision-hazard-decision
 information structure with hourly recourse:

\begin{subequations}
  \label{eq:WeeklyPlanningHourlyRecourseProblemFormulation}
  \begin{align}
    &
      \min_{ \va{\Controltilde}, \va{\RecourseControltilde}}
        \EE  \Biggl[
              \sum_{\week\in \WEEK}
              \widetilde{\InstantaneousCost}_{\week}\np{\va{\Stock}_{\whinf},
                                \va{\Controltilde}_{\closedWopen{\week}},
                                \va{\Uncertain}_{\openWclosed{\week}},
                                \va{\RecourseControltilde}_{\openWclosed{\week}}}
                  +
                  \FinalCost\np{\va{\Stock }_{\whlast}}			
                  \Biggr]
                  \\
      & \suchthat \eqsepv
      \forall \wh \in \timeline
      \nonumber
      \\
      &
      \va{\Stock}_{\winfhinf} = \va{\Uncertain}_{\winfhinf} \eqfinv
      \\
      &
      \va{\Stock}_{\wsuccessorhinf} =\widetilde{\Dynamics}_{\week}\np{\va{\Stock}_{\whinf}, \va{\RecourseControltilde}_{\openWclosed{\week}}}
              \eqfinv
      \\
      &
      \widetilde{\NodeBalance}_{\closedWopen{\week}}\bp{\va{\Controltilde}_{\closedWopen{\week}},
                                \va{\Uncertain}_{\openWclosed{\week}} ,  \va{\RecourseControltilde}_{\openWclosed{\week}}}
                                \geq
                                \va{0} \eqfinv \label{eq:balanceHourlyRecourse}
                                \\
      &\sigma\bp{\va{\Controltilde}_{\closedWopen{\week}} } \subset
      \sigma\bp{\va{\Uncertain}_{\winfhinf},
                  \va{\Uncertain}_{\openWclosed{\binf{\week}}}, \dots,
                  \va{\Uncertain}_{\openWclosed{\week\predecessor}} }
                  \eqfinv
      \\
      &\sigma\bp{\va{\RecourseControltilde}_{\wh\successor} } \subset
      \bsigmaf{
        \va{\Uncertain}_{\winfhinf},
        \va{\Uncertain}_{\openWclosed{\binf{\week}}},
        \dots,
        \va{\Uncertain}_{\openWclosed{\week\predecessor}},
        \va{\Uncertain}_{\whinfsuccessor},
        \dots,
        \va{\Uncertain}_{\wh\successor}
        }
            \eqfinp
   \end{align}
  \end{subequations}
  The final cost $\FinalCost\np{\va{\Stock }_{\whlast}}$ is used to give value to the energy in the storage at
  the end of the yearly period.

  \subsection{Bellman equations in decision-hazard-decision information structure with hourly recourse}
  \label{Bellman_equations_in_decision-hazard-decision_information_structure_with_hourly_recourse}

  Here, we provide Bellman equations corresponding to the
   stochastic multistage optimization
   problem~\eqref{eq:WeeklyPlanningHourlyRecourseProblemFormulation}.
 Their derivation is not given, as it is lengthy, but can be established from~\cite{Carpentier-Chancelier-DeLara-Martin-Rigaut:2023}.
  Defining the weekly state~$\statew=\stock_{\whinf}$, that is, the storage level at the beginning of the week,
  we write the weekly Bellman equations~\eqref{eq:BellmanHourlyRecourse}
  corresponding to~\eqref{eq:WeeklyPlanningHourlyRecourseProblemFormulation} as
\begin{subequations}
      \label{eq:BellmanHourlyRecourse}
  \begin{align}
    &\bBellman{\weeklast\successor}{\WPHR}{\state_{\weeklast\successor}} =  \FinalCost\np{\va{\Stock }_{\whlast}}	\\
    &\bBellman{\week}{\WPHR}{\state_{\week}} =
     \nonumber \\
       & \hspace{1cm}
            \min_{\control_{\closedWopen{\week}}}
            \hspace{0.1cm}     \EE  \Biggl[
            \nonumber  \\
              & \hspace{1.25cm}
               \min_{\recoursecontrol_{\whinfsuccessor}} \EE  \biggl[
                \nonumber   \\
                & \hspace{1.5cm}
                 \min_{\recoursecontrol_{\whinfsuccessor\successor}} \EE  \Bigl[
                  \nonumber  \\
               &
               \hspace{2.5cm}   \vdots
               \nonumber        \\
                     & \hspace{2.3cm}  \min_{\recoursecontrol_{\whsup}} \EE\bigl[
                      \nonumber       \\
                      &\hspace{2.6cm}
                                   \min_{\recoursecontrol_{\wsuccessorhinf}}
                         \bigl\{
                          \widetilde{\InstantaneousCost}_{\week}\np{{\statew},
                          {\controltilde}_{\closedWopen{\week}},
                          {\uncertain}_{\openWclosed{\week}},
                          {\recoursecontroltilde}_{\openWclosed{\week}}}
                 +  \nBellman{\week\successor}{\WPHR}{\widetilde{\Dynamics}_{\week}\np{\statew,
                                                                                                            \recoursecontroltilde_{\openWclosed{\week}}} }
                \bigr\}
                \nonumber    \\
                & \hspace{5.6cm} \mid\np{\uncertain_{\whinfsuccessor}, \uncertain_{\whinfsuccessor\successor} ,\dots, \uncertain_{\whsup} }
               \bigr]
               \nonumber    \\
               &
               \hspace{8.2cm}   \vdots
               \nonumber     \\
               & \hspace{7.5cm} \big| \np{\uncertain_{\whinfsuccessor}, \uncertain_{\whinfsuccessor\successor}} \Bigr]
               \nonumber    \\
               & \hspace{9.2cm} \Big|  \uncertain_{\whinfsuccessor} \biggr]
               \nonumber    \\
               &\hspace{10.5cm}  \Biggr] \eqfinv
            \label{eq:BellmanHourlyRecourseB}
\end{align}
\end{subequations}
where all the minimizations in Equation~\eqref{eq:BellmanHourlyRecourseB} are computed subject to the constraint~\eqref{eq:balanceHourlyRecourse}.

If the sequence $\bp{\va{\Uncertain}_{\openWclosed{\binf{\week}}}, \dots,
\va{\Uncertain}_{\openWclosed{\week}} , \dots, \va{\Uncertain}_{\openWclosed{\bsup{\week}}}}$
of uncertainties is weekly  independent,
the weekly Bellman equations provide an optimal solution for the problem~\eqref{eq:WeeklyPlanningHourlyRecourseProblemFormulation}.
We highlight that the hourly uncertainties
 \(\va{\Uncertain}_{\openWclosed{\week}}=  \bp{\va{\Uncertain}_{\whinfsuccessor}, \dots, \va{\Uncertain}_{\wh\successor}, \dots, \va{\Uncertain}_{\wsuccessorhinf}} \)
within the week  are  not assumed to be independent to get an optimal solution.
Under this independence assumption,
it is well known that, for all~$\week\in \WEEK$, the function~$\nBellmanf{\week}{\WPHR}$ satisfies
\begin{subequations}
 \begin{align}
 \nBellman{\week}{\WPHR}{\statew} & =
    \min_{\substack{\va{\Controltilde}_{\closedWopen{\week}},\dots, \va{\Controltilde}_{\closedWopen{\bsup{\week}}}   \\
                    \va{\RecourseControltilde}_{\openWclosed{\week}},\dots, \va{\RecourseControltilde}_{\openWclosed{\bsup{\week}}} }}
        \EE \Biggl[
           \sum_{\week'=\week}^{\bsup{\week}}
           \InstantaneousCost_{\week'}\np{\va{\Stock}_{\wprimehinf},
           \va{\Controltilde}_{\closedWopen{\week'}},\va{\Uncertain}_{\openWclosed{\week'}},\va{\RecourseControltilde}_{\openWclosed{\week'}}
                       }
               +
               \FinalCost\np{\va{\Stock }_{\whlast}}			
               \Biggr]
               \\
   & \suchthat \eqsepv
   \forall \week' \in \ic{\week,\dots,\bsup{\week}}
   \nonumber
   \\
   &
   \va{\Stock}_{\wprimehinf} =  \statew \eqfinv
   \\
   &
   \va{\Stock}_{\wprimesuccessorhinf} = \Dynamics_{\week'}\np{\va{\Stock}_{\wprimehinf},
                                                           \va{\RecourseControltilde}_{\openWclosed{\week'}}
                                                        }		
           \eqfinv
   \\
   &
   \NodeBalance_{\closedWopen{\week'}}\np{\va{\Controltilde}_{\closedWopen{\week'}},
   \va{\Uncertain}_{\openWclosed{\week'}}, \va{\RecourseControltilde}_{\openWclosed{\week'}}}
   \geq
        \va{0} \eqfinv \\
   &
   \sigma\np{ \va{\Controltilde}_{\closedWopen{\week'}}}
                       \subset \bsigmaf{
                                   \va{\Uncertain}_{\winfhinf},
                                   \va{\Uncertain}_{\openWclosed{\binf{\week}}},
                                   \dots,
                                   \va{\Uncertain}_{\openWclosed{\wprimepredecessor}}
                                   }
                 \eqfinv
                 \\
   &
   \sigma\bp{\va{\RecourseControltilde}_{\wprimeh\successor} } \subset
      \bsigmaf{
        \va{\Uncertain}_{\winfhinf},
        \va{\Uncertain}_{\openWclosed{\binf{\week}}},
        \dots,
        \va{\Uncertain}_{\openWclosed{\wprimepredecessor}},
        \va{\Uncertain}_{\wprimehinfsuccessor},
        \dots,
        \va{\Uncertain}_{\wprimeh\successor}
        }
            \eqfinp
\end{align}

\end{subequations}
so that  the  value $ \nBellman{\week}{\WPHR}{\statew}$ of Bellman function is interpreted as the future optimal cost
when at week $\week$ the storage level is $\statew$ and the weekly decision-hazard-decision
information structure with hourly recourse is considered.
More details can be found in \cite[Proposition~10]{Carpentier-Chancelier-DeLara-Martin-Rigaut:2023}.

\subsection{Theoretical comparison between Bellman functions in $\HD$, $\DHD$ and $\WPHR$}
\label{subsection:TheoreticalBellmanComparisonHourlyRecourse}
 When we compare theoretically the Bellman functions $\nBellmanf{\week}{\HD}$, $\nBellmanf{\week}{\DHD}$ and $ \nBellmanf{\week}{\WPHR}$ , given
 by the Bellman equations~\eqref{eq:HDweeklyBellman},~\eqref{eq:DHDweeklyBellman} and~\eqref{eq:BellmanHourlyRecourse},
 we observe that the weekly hazard-decision  approach is a relaxation of the weekly decision-hazard-decision  approach, and this last one
 is a relaxation of the weekly decision-hazard-decision with hourly recourse approach
  with respect to the information constraint.

  In other words, at each stage of the stochastic multistage optimization problem, the decision maker has less information when changing from
   the $\HD$ to the $\DHD$ approach. The same happens when changing from
   the $\DHD$ to the $\WPHR$ approach.  Therefore, we have the following inequalities for all
$\week\in\WEEK$:
\begin{align}
  \label{eq:BellmanRelation2}
  \nBellmanf{\week}{\HD} \leq \nBellmanf{\week}{\DHD} \leq \nBellmanf{\week}{\WPHR}  \eqfinp
\end{align}


\begin{thebibliography}{10}

\bibitem{Bellman:1957}
R.~E. Bellman.
\newblock {\em Dynamic Programming}.
\newblock Princeton University Press, Princeton, N.J., 1957.

\bibitem{Bertsekas:2000}
D.~P. Bertsekas.
\newblock {\em Dynamic Programming and Optimal Control}.
\newblock Athena Scientific, Belmont, Massachusetts, second edition, 2000.
\newblock Volumes 1 and 2.

\bibitem{bezanson2017julia}
J.~Bezanson, A.~Edelman, S.~Karpinski, and V.~B. Shah.
\newblock Julia: A fresh approach to numerical computing.
\newblock {\em SIAM review}, 59(1):65--98, 2017.

\bibitem{Carpentier-Chancelier-DeLara-Martin-Rigaut:2023}
P.~Carpentier, J.-P. Chancelier, M.~De~Lara, T.~Martin, and T.~Rigaut.
\newblock Time block decomposition of multistage stochastic optimization
  problems.
\newblock {\em Journal of Convex Analysis}, 30(2), 2023.

\bibitem{Dellacherie-Meyer:1975}
C.~Dellacherie and P.~A. Meyer.
\newblock {\em Probabilit\'es et potentiel}.
\newblock Hermann, Paris, 1975.

\bibitem{Antares}
M.~Doquet, R.~Gonzalez, S.~Lepy, E.~Momot, and F.~Verrier.
\newblock A new tool for adequacy reporting of electric systems: Antares.
\newblock CIGRE, 2008.
\newblock 42nd International Conference on Large High Voltage Electric Systems.

\bibitem{dowson2020policy}
O.~Dowson.
\newblock The policy graph decomposition of multistage stochastic programming
  problems.
\newblock {\em Networks}, 76(1):3--23, 2020.

\bibitem{LEE2001269}
J.~Lee and D.~Wilson.
\newblock Polyhedral methods for piecewise-linear functions {I:} the lambda
  method.
\newblock {\em Discrete Applied Mathematics}, 108(3):269--285, 2001.

\bibitem{Lubin2023}
M.~Lubin, O.~Dowson, J.~{Dias Garcia}, J.~Huchette, B.~Legat, and J.~P. Vielma.
\newblock {JuMP~1.0}: {R}ecent improvements to a modeling language for
  mathematical optimization.
\newblock {\em Mathematical Programming Computation}, 2023.

\bibitem{Pacaud-DeLara-Chancelier-Carpentier:2022}
F.~Pacaud, M.~{De Lara}, J.-P. Chancelier, and P.~Carpentier.
\newblock Distributed multistage optimization of large-scale microgrids under
  stochasticity.
\newblock {\em IEEE Transactions on Power Systems}, 37(1):204--211, 2022.

\bibitem{Pereira-Pinto:1991}
M.~V.~F. Pereira and L.~M. V.~G. Pinto.
\newblock Multi-stage stochastic optimization applied to energy planning.
\newblock {\em Math. Program.}, 52:359--375, Oct. 1991.

\bibitem{Steeger2014}
G.~Steeger, L.~Barroso, and S.~Rebennack.
\newblock Optimal bidding strategies for hydro-electric producers: A literature
  survey.
\newblock {\em Power Systems, IEEE Transactions on}, 29:1758--1766, 07 2014.

\bibitem{Street-Valladao-Lawson-Velloso:2020}
A.~Street, D.~Vallad{\~a}o, A.~Lawson, and A.~Velloso.
\newblock Assessing the cost of the hazard-decision simplification in
  multistage stochastic hydrothermal scheduling.
\newblock {\em Applied Energy}, 280:115939, 2020.

\bibitem{Valladao-Silva-Poggi:2019}
D.~Vallad{\~a}o, T.~Silva, and M.~Poggi.
\newblock Time-consistent risk-constrained dynamic portfolio optimization with
  transactional costs and time-dependent returns.
\newblock {\em Annals of Operations Research}, 282:379--405, 2019.

\bibitem{zou2019stochastic}
J.~Zou, S.~Ahmed, and X.~A. Sun.
\newblock Stochastic dual dynamic integer programming.
\newblock {\em Mathematical Programming}, 175:461--502, 2019.

\end{thebibliography}
\end{document}